\documentclass[leqno,11pt]{article}
\usepackage{amsmath}
\usepackage{amsfonts}
\usepackage{amssymb}
\usepackage{graphicx}

\usepackage{euscript}

\usepackage[latin1]{inputenc}
\usepackage{amsmath}
\usepackage{amsfonts}

 \evensidemargin\oddsidemargin
\renewcommand{\theequation}{\thesection.\arabic{equation}}
\newtheorem{theorem}{Theorem}[section]
\newtheorem{lemma}[theorem]{Lemma}
\newtheorem{proposition}[theorem]{Proposition}

\newtheorem{remark}[theorem]{Remark}

\newcommand{\eqnsection}{
\renewcommand{\theequation}{\thesection.\arabic{equation}}
    \makeatletter
    \csname  @addtoreset\endcsname{equation}{section}
    \makeatother}
\eqnsection

\newtheorem{prop}[theorem]{Proposition}

\def\r{{\mathbb R}}
\def\e{{\mathbb E}}
\def\p{{\mathbb P}}
\def\q{{\mathbb Q}}
\def\z{{\mathbb Z}}

\newcommand{\Ray}{{\mbox{\tt Ray}}}
\newcommand{\IGWR}{{\mbox{\tt IGWR}}}
\newcommand{\IGW}{{\mbox{\tt IGW}}}
\newcommand{\GW}{{\mbox{\tt GW}}}
\newcommand{\T}{{\mathcal{T}}}
\newcommand{\reals}{\mathbb{R}}

\def\squarebox#1{\hbox to #1{\hfill\vbox to #1{\vfill}}}
\newcommand{\qed}{\hspace*{\fill}
\vbox{\hrule\hbox{\vrule\squarebox{.667em}\vrule}\hrule}\smallskip}

\title{Einstein relation for biased random walk on Galton--Watson trees}\author{Gerard Ben Arous\thanks{Courant institute, New York University, 251 
Mercer St., New York, NY 10012, U.S.A. Email: benarous@cims.nyu.edu.} \and 
Yueyun Hu\thanks{D\'{e}partement de Math\'{e}matiques, LAGA, 
Universit\'{e} Paris 13,
99 Av. J-B Cl\'{e}ment, 93430 Villetaneuse, FRANCE.
Email: yueyun@math.univ-paris13.fr}
\and
Stefano Olla\thanks{CEREMADE,
Place du Mar\'{e}chal de Lattre de TASSIGNY, 
F-75775 Paris Cedex 16, FRANCE.  
\emph{and} 
 INRIA, Projet MICMAC, Ecole des
 Ponts, 6 \& 8 Av. Pascal, 77455 Marne-la-Vall\'ee Cedex 2, France.
Email: olla@ceremade.dauphine.fr .
The work of this author was partially supported by the
  European Advanced Grant {\em Macroscopic Laws and Dynamical Systems}
  (MALADY) (ERC AdG 246953) and by grant ANR-2010-BLAN-0108 (SHEPI).
}\and
Ofer Zeitouni\thanks{School of Mathematics, University of Minnesota,
206 Church St. SE, Minneapolis, MN 55455, USA \emph{and}
 Faculty of Mathematics,
Weizmann
Institute, POB 26, Rehovot 76100, Israel. Email: zeitouni@math.umn.edu.
The work of this author was partially
supported by NSF grant DMS-0804133 and by a grant from the Israel Science Foundation.}}

\date{June 19, 2011. Revised November 27, 2011}

\begin{document}

\maketitle
\abstract{We prove the Einstein relation, relating the velocity under a small
perturbation to the diffusivity in equilibrium, for certain biased random
walks on Galton--Watson trees. This provides the first example where 
the Einstein relation is proved for motion in random media with arbitrary
deep traps.}

\section{Introduction}
Let $\omega$ be a rooted Galton--Watson tree with offspring distribution
$\{p_k\}$, where $p_0=0$, $m=\sum k p_k>1$ and $\sum b^k p_k<\infty$ for
some $b>1$. For a vertex $v\in \omega$, let $|v|$ denote the
distance of $v$ from the root of $\omega$. Consider a
(continuous--time) nearest-neighbor random walk $\{Y_t^\alpha\}_{t\geq 0}$ 
on $\omega$, 
which when at a vertex $v$, jumps with 
rate $1$ toward each child of  $v$ and at rate $\lambda=
\lambda_\alpha = me^{-\alpha}$,
$\alpha\in \r$,
toward the parent of $v$.

It follows from \cite{RL} that if $\alpha=0$, the random walk
 $\{Y_t^\alpha\}_{t\geq 0}$ is, for 
almost every tree $\omega$, null recurrent (positive recurrent for $\alpha<0$,
transient for $\alpha>0$). Further, an easy adaptation 
of \cite{PZ} 
shows that $|Y_{[nt]}^0|/\sqrt{n}$ satisfies a (quenched, and hence also
annealed) invariance principle (i.e., converges weakly to a multiple of 
the  absolute value
of a Brownian motion),
with diffusivity 
\begin{equation}
\label{eq-of1}
{\cal D}^0=\frac{2m^2(m-1)}{\sum k^2 p_k-m}\,.
\end{equation}
(Compare with \cite[Corollary 1]{PZ}, and note that the factor $2$
is due to the speed up of the continuous--time walk relative 
to the discrete--time walk considered there. See \eqref{eq-diff0}
below and also the derivation in \cite{DS}.) 
On the other
hand, see \cite{LPP}, when $\alpha>0$, 
$|Y_t^\alpha|/t\to_{t\to\infty} \bar
v_\alpha>0$, almost surely, with $\bar v_\alpha$
deterministic. A consequence of
our main result, Theorem \ref{theo-main} below, is the following.
\begin{theorem}[Einstein relation]
\label{theo-GW}
With notation and assumptions as above, 
\begin{equation}
\label{ER-GW}
\lim_{\alpha\searrow 0} \frac{\bar 
v_\alpha}{\alpha}= \frac{ {\cal D}^0}{2}\,.
\end{equation}
\end{theorem}
The relation \eqref{ER-GW} is known as an {\it Einstein relation}.
It is straight forward to verify that for homogeneous random walks
on $\z_+$ (corresponding to deterministic Galton--Watson trees, that is, those
with $p_k=1$ for some $k\geq 1$), the Einstein relation holds.

In a weak limit (velocity rescaled with time) the Einstein relation is
proved in a very general setup by Lebowitz and Rost
(cf. \cite{lero}). See also \cite{BD} for general
fluctuation-dissipation relations. 

For the tagged particle in the symmetric exclusion process, the
Einstein relation has been proved by Loulakis in $d\ge 3$  \cite{LOU}.
The approach of  \cite{LOU}, based on perturbation theory and
transient estimates, was adapted for bond diffusion in $\mathbb
Z^d$ for special environment distributions (cf. \cite{OK}). For mixing
dynamical random 
environments with spectral gap, a full perturbation expansion can also be
proved (cf. \cite{OK1}).

For a diffusion in random potential, the recent \cite{GMP} proves the
 Einstein relation by following the strategy of 
 \cite{lero}, adding to it 
 a  good control (uniform in the environment) of suitably defined
 regeneration times in the transient 
 regime.  
A major difference in our setup is the possibility of having 
``traps'' of arbitrary strength 
in the environment; in particular, the presence of such
traps does not allow one to obtain estimates on regeneration times that
are uniform in the environment, and we 
have been unable to
obtain sharp enough estimates on regeneration times
that would allow us to mimic the strategy in
\cite{GMP}. On the other hand, the tree structure allows us to develop
some estimates directly for hitting times via recursions, see Section
\ref{sec:drift-towards-desc}. We emphasize that our work is (to the best of our 
knowledge) the first in which an Einstein relation is rigorously
proved for motion in random environments with arbitrary strong traps.

In order to explore the full range of parameters $\alpha$, we will work in
a more general context than that described above, following 
\cite{PZ}.  This is described next.

Consider infinite trees $\T$ with no leaves, equipped 
with one (semi)-infinite
directed path, denoted $\Ray$, starting  from a distinguished 
vertex called the {\bf root} and denoted $o$. We call such a tree a
{\it marked tree}.
Using $\Ray$, 
we define in a natural way the offsprings of a vertex $v\in \T$,
and denote by $D_n(v)$ the collection of vertices that are descendants of $v$ 
at distance $n$ from $v$, with $Z_n(v)=|D_n(v)|$. See \cite[Section 4]{PZ}
for precise definitions.
For any vertex $v\in {\T}$,
we let $d_v$ denote the number of offspring of $v$, and write
${\buildrel
\leftarrow \over v}$ for the parent of $v$.  
Finally, we write $\rho(v)$ for the 
{\it horocycle} distance of $v$ from the root $o$. Note that $\rho(v)$ is 
positive if $v$ is a a descendant of $o$ and negative if it is an ancestor
of $o$.

Let $\Omega_T$ denote the space of marked trees. As in \cite{PZ} and motivated
by \cite{LPP}, given an offspring distribution $\{p_k\}_{k\geq 0}$ satisfying
our general assumptions,
we introduce a  reference
probability $\IGW$ on $\Omega_T$, as follows.
Fix   the root  $o$  and a semi-infinite ray, denoted $\Ray$,
emanating from it. Each vertex $v\in \Ray$ with $v\neq o$
is assigned 
independently a size-biased number of offspring, that
is $P_{\IGW}(d_v=k)=kp_k/m$, one of which is
identified with the descendant of $v$ on $\Ray$.
To each offspring of $v\neq o$ not on $\Ray$,
and to $o$,
one attaches an independent Galton-Watson tree of offspring
distribution $\{p_k\}_{k\geq 0}$.
Note that $\IGW$ makes
the collection $\{d_v\}_{v\in \cal T}$ independent.
We denote expectations with respect to
$\IGW$ by $\langle\cdot\rangle_0$ (the reason for the notation
will become apparent in Section \ref{sec-env} below).

As mentioned above and in
contrast with \cite{LPP} and \cite{PZ},
it will be convenient to work in continuous time, because
it slightly simplifies the formulas
(the
adaptation needed to transfer the results to 
the discrete time setup of \cite{LPP} are straight-forward). 
For background,
we refer to \cite{DS}, where the results in \cite{LPP} and \cite{PZ}
are transferred to continuous time, in the more general setup of multi-type
Galton--Watson trees. Given a marked tree
${\T}$ and $\alpha\in \reals$, we define 
an $\alpha$-biased random walk $\{X_t^\alpha\}_{t\geq 0}$ 
on $\T$ as the continuous 
time Markov process with state space the vertices of ${\T}$, 
$X_0^\alpha=o$, and 
so that when at $v$,
the jump rate is $1$ toward   each of the descendants of $v$,
and the jump rate is $e^{-\alpha} m$ toward
the parent ${\buildrel
\leftarrow \over v}$. More explicitly, the generator of
the random walk $\{X_t^\alpha\}_{t\geq 0}$
can be written as
\begin{equation}
  \label{eq:2}
  \mathcal L_{\alpha, {\T}} F(v) = \sum_{x\in D_1(v)} 
  \left( F(x) - F(v) \right) + e^{-\alpha} m \left(F({\buildrel
\leftarrow \over v}) -
    F(v)\right)\,. 
\end{equation}
Alternatively,
\begin{equation}
  \label{eq:2bis}
  \mathcal L_{\alpha,{\T}} F(v) = - m\partial_{-}^* \partial_{-} F(v) +
  (e^{-\alpha} -1) m \partial_- F(v)
\end{equation}
where $\partial_- F(v) = F({\buildrel
\leftarrow \over v}) - F(v)$ and
\begin{equation}
  \label{eq:3}
  \partial_{-}^* F(v) = \frac 1m \sum_{x\in D_1(v)} F(x) - F(v)
\end{equation}
Note that if
$\alpha < 0$ the (average) drift is towards the ancestors, 
whereas if $\alpha > 0$ the (average) drift is towards the children.
As in \cite{LPP} and \cite{PZ}, we have that
\begin{equation}
\lim_{t\to\infty} \frac{\rho(X_t^\alpha)}{t}\to_{t\to \infty} 
v_\alpha\,, \quad {\rm \IGW-a.s.}
\end{equation}
It is easy to verify that when 
$\alpha>0$, then $v_\alpha=\bar v_\alpha$, and that $\mbox{\rm 
sign}(v_\alpha)=
\mbox{\rm sign}(\alpha)$.
Further, we have, again from 
\cite{PZ}, that $\rho(X_{[nt]}^0)/\sqrt{n}$ satisfies
the invariance principle (that is, converges weakly
to a Brownian motion), with diffusivity constant ${\cal D}^0$ as in
\eqref{eq-of1}.

Our main result concerning walks on IGW-trees is the following.
\begin{theorem}
\label{theo-main}
With assumptions as above,
\begin{equation}
\label{eq-of3}
\lim_{\alpha\to 0} \frac{v_\alpha}{\alpha} = \frac{ {\cal D}^0}{2}\,.
\end{equation}
\end{theorem}

\begin{remark}
	It is natural to expect that the Einstein relation holds in 
	many related models, including Galton--Watson trees with only
	moment bounds on the offspring distribution, multi-type Galton--Watson
	trees as in \cite{DS}, and walks in random environments on 
	Galton--Watson trees, at least in the regime
	where a CLT with non-zero variance holds, see \cite{Far}.
	We do not explore these extensions here.
\end{remark}

The structure of the paper is as follows. In the next section, we consider
the case of $\alpha<0$, exhibit an invariant measure for the environment
viewed from the point of view of the particle, and use it to prove 
the Einstein relation when $\alpha\nearrow 0$. Section
\ref{sec:drift-towards-desc} deals with the harder case of $\alpha\searrow 0$.
We first prove an easier Einstein relation (or linear response) concerning
escape probabilities of the walk, exploiting the tree structure to introduce 
certain recursions. Using that, we relate the Einstein relation for velocities 
to estimates on hitting times. A crucial role in obtaining
these estimates, and an alternative formula for
the velocity (Theorem \ref{theo-3.4}), is
obtained by 
the introduction, after \cite{FHS}, of a spine random walk, 
see Lemma \ref{L:zjbis}.

\section{The environment process, and proof of Theorem
\ref{theo-main} for $\alpha\nearrow 0$.}
\label{sec-2}
As is often the case when motion in random media is concerned, it is 
advantageous to consider the evolution, in $\Omega_T$,
of the environment from the point of view
of the particle. One of the reasons for our opting to work in
continuous time is that when $\alpha=0$, the invariant measure
for that (Markov) process is simply \IGW, in contrast with the more complicated
measure \IGWR  \ of \cite{PZ}. We will see that when $\alpha<0$, 
an explicit invariant measure for the environment viewed from
the point of view of the particle exists, and is absolutely continuous
with respect to \IGW.
\subsection{The environment process}
\label{sec-env}
For a given tree $\T$ and
$x\in {\T}$, let $\tau_x $ denote 
the shift that moves the root of $\T$ to
$x$, with \Ray\ shifted to start at $x$ in the unique way so
 that it differs from
\Ray\ before the shift by only finitely many vertices.  
Then $\tau_x \T$ is rooted at $x$ and has the same (nonoriented) edges as $\T$.  (A special role will be played by $\tau_x $ for $x\in D_1(o)$, and
by $\tau_x$ with $x={\buildrel
\leftarrow \over o}$. We use
$\tau^{-1} {\T}=\tau_{{\buildrel
\leftarrow \over o}} {\T}$ in the latter case.)
The environment process $\{ {\T}_t\}_{t\geq 0}$
is defined by ${\T}_t=\tau_{X_t} {\T}$.
It is straightforward to check that the environment
process is a Markov process. In fact, introducing  the operators
$$
Df({\cal T}) = f(\tau^{-1}\T) - f(\T)\,,
$$
we have that the adjoint operator (with respect to $\IGW$) is
\begin{equation*}
  D^* f(\T) = \frac 1m \sum_{x\in D_1(o)} f(\tau_x \T) - f(\T)
\end{equation*}
since 
\begin{equation*}
  \left< g Df \right>_0 =  \left< f D^* g \right>_0.
\end{equation*}
Notice that $D^* 1= d_o/m -1$. 

Define $W(v,n) = Z_n(v)/m^n$. Then $W(v,n)$ is a positive 
martingale that converges  to a random
variable denoted $W_v$.
Using the recursions
\begin{equation*}
  m W_v = \sum_{x\in D_1(v)} W_x, \qquad  
   m W(v,n) = \sum_{x\in D_1(v)} W(n-1,x),
\end{equation*}
we see that
$\langle W_v\rangle_0=1$ for $v\not\in\Ray$.
To simplify notation, we write
$W_{-j} = W_{v_j}$ with $v_j\in \Ray$ denoting the
$j$-th ancestor of $o$.
Since $W_o(\tau_x\T) = W_x(\T)$, we have that $D^* W_o = 0$.

The generator of the environment process is 
\begin{equation}
  \begin{split}
    L_\alpha f(\T)= \sum_{x\in D_1(o)}\left[ f(\tau_x \T) -
      f(\T) \right] + e^{-\alpha} m \left[ f(\tau^{-1} \T) -
      f(\T) \right] \\
    = - m D^* D f (\T) + (e^{-\alpha} - 1) m D f (\T)
    \label{eq:1}
  \end{split}
\end{equation}
The adjoint operator (with respect to \IGW) 
is $L_\alpha^* = -m D^* D +  (e^{-\alpha} - 1) m D^*$.
For any $\alpha\in \r$,
let $\mu_\alpha$ 
denote any stationary probability measure for $L_\alpha$, 
that is $\mu_\alpha$ satisfies, for any bounded measurable $f$,
\begin{equation*}
  \left< L_\alpha f \right>_\alpha = 0, 
\end{equation*}
where $\left< g\right>_\alpha=\int g d\mu_\alpha$.

Note that 
$\IGW$
is stationary and \emph{reversible} for $L_0$. Further, it
is ergodic for the environment process.
This is elementary to prove, since for any bounded function $f(\T)$
such that $L_0 f = 0$, we have that $\langle |Df|^2\rangle_0 = 0$, i.e. $f$ is
translation invariant for a.e. $\T$ with respect to $\IGW$,
i.e. constant a.e. .
Thus, necessarily, 
$\mu_0=\IGW$,
justifying
our notation $\langle \cdot\rangle_0=\langle \cdot\rangle_\IGW$.

In our setup, due to the existence 
of regeneration times for $\alpha\neq 0$ with bounded expectation, a general
ergodic argument ensures the existence of a stationary
measure $\mu_\alpha$, which however may fail in general 
to be absolutely continuous
with respect to $\IGW$, see \cite{LPP}.
Further, because $\IGW$ is ergodic and the random walk is elliptic, 
there is at most
one  $\mu_\alpha$ which is absolutely
continuous with respect to \IGW, since under any such 
$\mu_\alpha$, the process $\T_t^\alpha$ must be ergodic,
see e.g. \cite[Corollary 2.1.25]{stflour}
for a similar argument. As we now show, when $\alpha<0$, 
this stationary measure $\mu_\alpha$ with density with respect to \IGW\
can be constructed explicitly.

\begin{lemma}
	\label{lem-stat}
  For $\alpha < 0$, 
  the probability measure $\mu_\alpha = \psi_\alpha \mu_{0}$  
  where
  \begin{equation}
    \label{eq:7}
    \psi_\alpha(\T) = C_\alpha^{-1} Z_\alpha\, ,
  \end{equation}
\begin{equation}
	\label{C-comp}
        \begin{split}
          Z_\alpha &=\sum_{j=0}^\infty e^{j\alpha} W_{-j}(\T), \\
          C_\alpha &= \frac{(1 - b)e^{\alpha}m^{-1}}{1 - e^{\alpha}
            m^{-1}} + \frac{ b e^{\alpha}}{1 - e^{\alpha}} +1, \qquad
          b = \frac{\sum_k k^2 p_k - m}{m(m-1)}\,.
        \end{split}
\end{equation}
is stationary for $L_\alpha$. Furthermore
\begin{equation}
  \label{eq:12}
   \lim_{\alpha \nearrow 0} \psi_\alpha(\T) = 1 \qquad \mu_0-a.e. 
\end{equation}
\end{lemma}

\noindent
\textbf{Proof of Lemma \ref{lem-stat}:}
We show first that $C_\alpha$ provides the correct normalization. 
In fact, from the relation
\begin{equation}
	\label{eq-ws}
	W_{-j}=m^{-1}W_{-j+1}+m^{-1}\sum_{s\in D_1({v_{-j}}),
	s\not\in \Ray} W_s=:m^{-1}(W_{-j+1}+L_j)
\end{equation}
and since $\langle W_s\rangle_0 = 1$ if $s\not\in\Ray$, we obtain
\begin{equation}
  \label{eq:17}
  \langle W_{-j}\rangle_0 = m^{-1}\langle W_{-j+1}\rangle_0 + m^{-1}b (m-1), 
  \quad j\geq 1\,.
\end{equation}
Since $\langle W_o\rangle_0= 1$, we deduce that 
\begin{equation}
  \label{eq:20}
  \langle W_{-j}\rangle_0 = (1 - b) m^{-j} + b, \quad j\geq 0\,. 
\end{equation}
Thus,
\begin{eqnarray*}
	\langle \sum_{j=0}^\infty e^{j\alpha} W_{-j}\rangle_0&=&
	1+ b \sum_{j=1}^\infty e^{j\alpha}+(1-b)\sum_{j=1}^\infty
	e^{j \alpha }m^{-j}\\
	&=& 1+\frac{b e^{\alpha}}{1-e^{\alpha}}+
	\frac{(1-b)e^{\alpha}m^{-1}}{1-e^{\alpha}m^{-1}}\\
        &=& \frac{b}{1-e^{\alpha}}+
	\frac{1-b}{1-e^{\alpha}m^{-1}}=C_\alpha\,,
\end{eqnarray*}
as needed.

Note that the terms $L_j$ appearing in the right side
of (\ref{eq-ws})
are i.i.d..
Substituting and iterating,
we get
$$
W_{-k}=\frac{W_o}{m^k}+\frac{L_1}{m^k}+\frac{L_2}{m^{k-1}}+\cdots+
\frac{L_k}{m}\,.
$$ 
Therefore,
\begin{equation}
	\label{eq-ws1}
	\left(1-\frac{e^{\alpha}}{m}\right)Z_\alpha=
W_o
	+\frac{1}{m}
	\sum_{j=1}^\infty e^{\alpha j} L_j
	=:
	W_o+M_\alpha
	\,.
\end{equation}
Note that 
$M_\alpha$ is a weighted sum of i.i.d. random variables. Further, because
$\langle |D_1({v_{-j}})|\rangle_0=\sum k^2 p_k/m$ and 
$\langle W_s\rangle_0=1$, we have that
$\lim_{\alpha\nearrow 0}
\vert \alpha\vert \, \langle M_\alpha\rangle_0=(\sum k^2 p_k-m)/m^2:=\bar C$, and that
$\mbox{\rm Var}_{\IGW}(M_\alpha)=O({1\over \vert \alpha\vert})$. It then follows 
(by an interpolation argument) that 
that
$$\lim_{\alpha\nearrow 0}  \vert \alpha\vert \, M_\alpha =\bar C\,, \quad \IGW-a.s.$$
Substituting in (\ref{eq-ws1}), this yields 
\begin{equation*}
  \lim_{\alpha\nearrow 0} \alpha Z_\alpha = b \quad \IGW-a.s.
\end{equation*}
and \eqref{eq:12} follows.

We next verify that $L_\alpha^* \psi_\alpha=0$.
Since $W_{-j}(\tau^{-1}\T) = W_{-j-1}(\T)$, 
we have
\begin{eqnarray*}
    D\psi_\alpha &=&  C_\alpha^{-1} \sum_{j=0}^\infty e^{j\alpha}
    \left(W_{-j-1}(\T) - W_{-j}(\T) \right) \\
    &=& C_\alpha^{-1} \left(\sum_{j=1}^\infty e^{(j-1)\alpha}
    W_{-j}(\T)  - \sum_{j=0}^\infty e^{j\alpha} W_{-j}(\T) \right)\\
    &=&   C_\alpha^{-1} \sum_{j=0}^\infty  
    (e^{(j-1)\alpha} - e^{j\alpha})
     W_{-j}(\T) -  C_\alpha^{-1} e^{-\alpha} W_o\\
     &=&  (e^{-\alpha} -1)  C_\alpha^{-1} 
     \sum_{j=0}^\infty  e^{j\alpha}
     W_{-j}(\T) - C_\alpha^{-1} e^{-\alpha} W_o \\
     &=&  (e^{-\alpha} -1) \psi_\alpha -  C_\alpha^{-1} e^{-\alpha} W_o\,.
\end{eqnarray*}
Since $D^* W_o = 0$, we have
\begin{equation*}
  D^* \left( D\psi_\alpha - (e^{-\alpha} -1) \psi_\alpha \right) = 0
\,,\end{equation*}
i.e.
\begin{equation*}
  L_\alpha^* \psi_\alpha =  0.
\end{equation*}
\qed

We can now provide the proof of Theorem \ref{theo-main}
in case $\alpha\nearrow 0$.\\
\textbf{ Proof of Theorem \ref{theo-main} when $\alpha\nearrow 0$:}
We begin with the computation of $v_\alpha$. Because $\mu_\alpha$ is
ergodic and absolutely continuous with respect to
\IGW, we have that $v_\alpha$ equals the average drift (under 
$\mu_\alpha$) at $o$, that is
\begin{eqnarray*}
    v_\alpha &=& m \langle  \frac{d_o}m -e^{-\alpha} \rangle_\alpha = 
    m \langle  D^*1 \rangle_\alpha -m (e^{-\alpha} -1) =
     m \langle  D^*1 \psi_\alpha \rangle_0 -m (e^{-\alpha} -1)\\
     &=&  
     m \langle D\psi_\alpha \rangle_0 -m (e^{-\alpha} -1)
     \\
     &=&  m \langle (e^{-\alpha} -1) \psi_\alpha -
     C_\alpha^{-1} e^{-\alpha}   W_o\rangle_0 - m (e^{-\alpha} -1)
     = - m  C_\alpha^{-1} e^{-\alpha}.   
\end{eqnarray*}
 
Thus,
 \begin{equation}
	 \label{eq-vellim}
	 \lim_{\alpha\nearrow 0} \frac{v_\alpha}{|\alpha|}=
  - \frac{m^2 (m-1)}{\sum k^2 p_k-m}\,.
  \end{equation}

  It remains to compute the diffusivity ${\cal D}^0$
  when $\alpha=0$. Toward this end, one simply repeats 
the 
computation in \cite[Corollary 1]{PZ}.
One obtains
that the diffusivity is 
\begin{equation}
	\label{eq-diff0}
	{\cal D}^0=\frac{\langle mW_o^2+\sum_{s\in D_1(o)} W_{s}^2\rangle_0}{
\langle W_o^2\rangle_0^2}\,.
\end{equation}
From the definitions we have that
$\langle W_o^2
\rangle_0= (\sum k^2p_k-m)/m(m-1)$ (see \cite[(2)]{PZ}), and thus
\begin{equation}
	\label{eq-diffusivity}
{\cal D}^0=\frac{2m^2(m-1)}{\sum k^2 p_k-m}\,.
\end{equation}
Together with \eqref{eq-vellim},
this completes the proof of Theorem \ref{theo-main} when 
$\alpha\nearrow 0$. \qed

\begin{remark}
	Note that the construction above fails for $\alpha>0$, because then
	$Z_\alpha$ is not defined. The case $\alpha=\infty$ is however
	special. In that case,
the generator is 
\begin{equation}
  \begin{split}
    L_{\infty} f(\T)= \sum_{x\in D_1(o)}\left[ f(\tau_x \T) -
      f(\T) \right]
    \label{eq:-infty}
  \end{split}
\end{equation}
In particular, one can  verify that the measure defined by
  $d\mu_{\infty} /  d\mu_{GW}=1/(Cd_o)$ with 
$C = \sum_k k^{-1}p_k$ and $\mu_{GW}$  the ordinary Galton--Watson
	measure \GW \ (defined as \IGW \ but with the standard Galton--Watson
	measure also for vertices on \Ray), is a stationary measure, and that
$v_\infty =\langle d_o\rangle_\infty=C$. It follows that
the natural invariant measure is not absolutely continuous
with respect to $\IGW$. 
\end{remark}
\begin{remark}
	For $\alpha<0$ one can construct
	other invariant measures, that of course are singular with
	respect to \IGW. A particular family of such measures 
	is absolutely continuous with respect to the ordinary Galton--Watson
	measure \GW 
        . Indeed, one can verify that the positive function
	\begin{equation}
  \label{eq:33}
  \psi(\T) = C \sum_{j=1}^\infty \prod_{i=1}^{j-1}  
  \frac{d_{-i}}{me^{-\alpha}} 
  = C \sum_{j=1}^\infty (me^{-\alpha})^{-j+1} \prod_{i=1}^{j-1} d_{-i}  \,,
\end{equation}
with $C=1-e^{\alpha}$, satisfies $\int \psi d\GW=1$ and $\psi d\GW$
is an invariant
measure for $L_\alpha$. One can also check that the Einstein relation
\eqref{ER-GW} is not satisfied under this measure, emphasizing the
role that the measure \IGW \ plays in our setup.
\end{remark}

\section{Drift towards descendants: proof of Theorem \ref{theo-main}
  for $\alpha \searrow 0$}
\label{sec:drift-towards-desc}

In the case $\alpha >0$ we cannot find an explicit expression for the
stationary measure so we have to proceed in a different way.
We first prove another form of
the Einstein relation 
in terms of
the \emph{escape probabilities} (probability of never returning to the origin).

Because we consider 
the case $\alpha>0$,
 there
is no difference between considering the walk under the Galton--Watson
tree or under \IGW \ -- \ the limiting velocity is the same, i.e.
$v_\alpha=\bar v_\alpha$. Thus, we only consider
the walk $\{Y_t^\alpha\}_{t\geq 0}$ below.

Our approach is to provide an alternative formula for 
the speed $v_\alpha$, see Theorem \ref{theo-3.4} below, which is valid
for all $\alpha>0$ small enough. In doing so, we 
will take advantage of certain recursions, and of the 
{\it spine random walk} associated with the walk on the Galton--Watson tree,
see Lemma \ref{L:zjbis}.

We recall our standing assumptions:
$p_0=0$, $m>1$, and $\sum b^k p_k<\infty$ for some $b>1$.
We will throughout drop the superscript $\alpha$ from the notation
when it is clear
from the context, writing e.g. $Y_t$ for $Y_t^\alpha$.
To introduce our recursions,
define $ T(x):= \inf\{t\ge 0: Y_t=x\}$ and $\tau_n:=
 \inf\{t\ge 0: |Y_t|=n\}$.
 For a given tree $\omega$,
 we write $P_{x,\omega}$ for the law of $Y_t$ with $Y_0=x$.
 For $0 < |x| \le n$, define \begin{eqnarray*}
 \beta_n(x) &:=& P_{x, \omega} \Big( T({\buildrel
\leftarrow \over x}) > \tau_n\Big), \quad
 \beta(x)  :=  P_{x, \omega} \Big( T({\buildrel
\leftarrow \over x})= \infty \Big)  , \\
\gamma_n(x)&:=& E_{x, \omega} \left( \tau_n \wedge T({\buildrel
\leftarrow \over x}) \right).
 \end{eqnarray*}

We study the recursions for $\beta_n$ and $\gamma_n$.
By the Markov property of $P_{x,\omega}$, 
 for $|x|<n$, \begin{eqnarray*} \gamma_n(x)&=&
\frac{1}{d_x+\lambda}+
\sum_{i=1}^{d_x} {1\over \lambda + d_x} E_{x_i, \omega} \left(
\tau_n \wedge T({\buildrel \leftarrow \over x}) \right) \\
&=&
\frac{1}{d_x+\lambda}+
 \sum_{i=1}^{d_x} {1\over \lambda + d_x} \left(  E_{x_i,
\omega} \left( \tau_n \wedge T(x) \right) +  P_{x_i, \omega}( T(x) <
\tau_n) \,  E_{x, \omega} \left( \tau_n \wedge T({\buildrel
\leftarrow \over x}) \right) \right),
\end{eqnarray*}

\noindent which implies that $$ \gamma_n(x)=
\frac{1}{d_x+\lambda}+
 \sum_{i=1}^{d_x}
{1\over \lambda + d_x} \left( \gamma_n(x_i) + (1- \beta_n(x_i))
\gamma_n(x)\right).$$ Hence for any $0<|x|< n$,  $$ \gamma_n(x)= {
1 + \sum_{i=1}^{d_x} \gamma_n(x_i) \over \lambda +
\sum_{i=1}^{d_x} \beta_n(x_i)},$$

\noindent with boundary condition $\gamma_n(x)= 0$ for any $|x|=n$. We take   the above equality as the definition of  $\gamma_n(o)$.   Similarly, we have
$$ \beta_n(x)= { \sum_{i=1}^{d_x} \beta_n(x_i) \over
\lambda + \sum_{i=1}^{d_x} \beta_n(x_i)} , \qquad 0< \vert x \vert < n, $$ with $\beta_n(x)=1$ if
$|x|=n$, and we define $\beta_n(o)$  so that the above equality holds for $x=o$. Finally, we let $\beta(o)= \lim_{n \to \infty} \beta_n(o)$ (the limit of the monotone sequence $\beta_n(o)$).

\begin{proposition}\label{eresc}
As $\alpha \searrow
0$,  $\alpha^{-1} \beta(o)$ converges in law and in expectation to a
  random variable $Y$ 
  such that
  \begin{equation}
    \label{eq:4}
    \mathbb E(Y) = \frac{m(m-1)}{\mathbb E(d_o^2 - d_o)}  = 
    \frac{\mathcal D^0}{2m}    
  \end{equation}
\end{proposition}
This is a form of Einstein relation, as linear response for the escape
probability. 
The law of $Y$ can be identified, see the end of the proof of Proposition \ref{eresc}.

{\noindent \bf Proof:} We clearly have that with 
 $  B(x):= { 1\over \lambda} \sum_{i=1}^{d_x}  \beta(x_i)$, it holds that
 $$ \beta(x)= 
 { B(x) \over 1+ B(x)}, \qquad \forall x \not= o,$$
and
\begin{equation} \label{bx} B(x)= {1\over \lambda} \sum_{ i=1}^{d_x} { B(x_i) \over 1+ B(x_i) }, \qquad \forall x\in {\cal T}
	. \end{equation}

 \noindent Notice that all $B(x)$ are distributed as 
 some random variable, say $B$, and conditionally on $d_x$ and on the tree 
 up to generation $\vert x \vert$, the variables
 $B(x_i), 1\le i \le d_x$ are i.i.d. and distributed as $B$.  
 It follows that  
 \begin{equation} \label{eb1}\e (B)= e^\alpha 
	 \e { B\over 1+ B}, \end{equation} 
	 and 
         \begin{equation}
 \e ( B^2) = { 1\over \lambda^2} 
	 \left( m \e \left( ( { B\over 1+B})^2 \right) + 
	 \e (d_o(d_o-1))\,    \left( \e ( { B\over 1+B})  \right)
         ^2\right).\label{eq:5} 
     \end{equation}

 \noindent  For any nonnegative r.v. $Z \in L^2$, let us denote by  $\{Z\}:= {Z \over \e(Z)}$. By concavity,   the  following   inequality holds (see e.g. \cite{PP}, Lemma 6.4):
 $$  \e \left( \{ { Z \over 1+ Z} \}^2\right) \le \e \left( \{ 
 Z\}^2\right) . $$

\noindent  By (\ref{eb1}) and (\ref{eq:5}), we get that 
 $$  \e \left( \{ B\}^2 \right) = { 1\over m} \e 
 \left( \{ { B \over 1+ B} \}^2\right) +  { \e(d_o(d_o-1))\over m^2} 
 \le { 1\over m} \e \left( \{ { B } \}^2\right) +  { \e(d_o(d_o-1))
 \over m^2},$$

 \noindent which yields that the second moment  of $B$ is uniformly
 bounded by the square of $\e(B)$: for any $0<\alpha$, $$ \e(B^2) 
 \le { \e(d_o(d_o-1)) \over m -1} \, ( \e(B))^2.$$

 \noindent
 By (\ref{eb1}),  $ e^{-\alpha} \e(B)= \e (B)- \e ({B^2\over 1+B})$, hence $  (1-e^{-\alpha}) \e(B)= \e ( { B^2 \over 1+ B}) 
 \le \e(B^2) \le { \e(d_o(d_o-1)) \over m -1} \, 
 ( \e(B))^2.$ It follows that  $$ \e(B) \ge  
 { m(m-1) \over \e(d_o(d_o-1))
 }\, 
 (1-e^{-\alpha}).$$ On the other hand, by Jensen's inequality, 
 $  \e (B) = e^\alpha \e { B\over 1+ B} 
 \le e^{\alpha}   { \e(B)\over 1+ \e(B)}, $  
 which implies that  $$ \e(B) \le (e^\alpha-1).$$
 Therefore, ${ B/ \alpha} $ is tight as $\alpha \searrow
 0$. In particular, 
 for some sub-sequence $\alpha\searrow
 0$, ${ B(x)/ \alpha} $ converges in law to some $Y(x)$. 
 Since ${B / \alpha}$ is bounded in $L^2$ uniformly 
 in $\alpha>0$ in a neighborhood of $0$,
 we deduce  from  (\ref{bx}) that 
 \begin{equation}
	 \label{late-night}
	 Y {\buildrel \mbox{d}  \over =} {1\over m} 
 \sum_{i=1}^N  Y_i ,
 \end{equation}
 where $N$ is distributed
 like $d_o$ and,
 conditionally on $N$, 
 $(Y_i)$ are i.i.d and distributed as $Y$; 
 moreover $ \e(Y)= \lim_{\alpha\searrow 0} 
 \e({B\over \alpha}) >0$ (the limit along the same sub-sequence).  

 Dividing
$$  
(1-e^{-\alpha}) \e(B)=  \e ( { B^2 \over 1+ B}) \le \e(B^2) 
$$ 
by $\alpha^2$, we get that $ \e (Y)= \e(Y^2) $. 
The same operation in \eqref{eq:5} gives
\begin{equation}
  \label{eq:6}
  E(Y)^2 = \frac{m(m-1)}{\mathbb E(d_o (d_o -1))} E(Y^2).
\end{equation}
Putting these together we obtain  $E(Y)= {{ \cal D}^o\over 2 m}$. 

On the other hand,  it is known (see e.g. \cite[Theorem 16]{AB})  that the law of $Y$ satisfying
	\eqref{late-night} is determined up to a multiplicative constant, and 
	therefore $Y$ equals in distribution $a W_o$ for some constant $a$. The 
	equality $EY=EY^2$ then implies that $Y$ equals in distribution
	$W_o/E(W_o^2)$.  Since all possible limits in law are the same, 
we  get that   
$ \beta(o)/\alpha$ converges in law to $W_o/E(W_o^2)$. $\Box$


\bigskip
\noindent We return to the proof of the Einstein relation concerning
velocities. Recall that a {\it level regeneration time}
is a time for which the random walk hits a fresh level and never 
backtracks, see e.g. \cite{DGPZ} for the definition and basic properties.
(Level regeneration times are related to, but different from,
the regeneration times introduced in \cite{LPP}.) In particular, see
\cite[Section 4]{DGPZ} and \cite[Section 7]{PZ}, the differences
of adjacent 
regeneration times form an i.i.d. sequence,
with all moments bounded.
Since $\gamma_n(x)$ is smaller than the $n$-th level regeneration
time (started at $x$), it follows that the sequence 
$\gamma_n(x)/n$ is uniformly integrable (under the measure
$\GW\times P_{x,\omega}$), and therefore,
the convergence in the forthcoming  \eqref{eq-gammalim} holds also in expectation: \begin{equation}\label{eq-gamma2} \lim_{n\to\infty} {\e [\gamma_n (o)] \over n} = { \e( \beta(o) )\over v_\alpha}.
\end{equation}

Since $  \gamma_n(x)= E_{x, \omega} \left( \tau_n 1_{(
\tau_n < T({\buildrel \leftarrow \over x}))} \right) + O(1)$
and  ${\tau_n \over n} \to {1\over v_\alpha},$ $P_{x, \omega}$
a.s. and in  $L^1$ (the latter follows at once from the integrability
of regeneration times mentioned above, as $\tau_n$ is bounded above
by the $n$th regeneration time),
we get that for $x$ fixed, 
\begin{equation}
	\label{eq-gammalim}
	{ \gamma_n (x)
\over
n} \, \to_{n\to\infty} \, {1\over v_\alpha} P_{x, \omega}\Big( T({\buildrel
\leftarrow \over x}) = \infty\Big)= {\beta(x) \over v_\alpha},
\qquad GW\, a.s.
\end{equation}

So all we need to prove in order to have the Einstein relation for velocities,
is that 
\begin{equation}
  \label{eq:9}
  \lim_{\alpha\to 0} \lim_{n\to \infty} \frac 1n \mathbb E
  \left(\gamma_n(o)\right) = \frac 1m .
\end{equation}
 To this end, define $$ B_n(x):= {1\over \lambda} \sum_{i=1}^{d_x}
\beta_n(x_i), \qquad \Gamma_n(x):=  \sum_{i=1}^{d_x}
\gamma_n(x_i), \qquad |x|<n.$$

Note that showing \eqref{eq:9} is equivalent to proving that
\begin{equation}
  \label{eq:9a}
  \lim_{\alpha\to 0} \lim_{n\to \infty} \frac 1n \mathbb E
  \left(\Gamma_n(o)\right) = 1 .
\end{equation}

For $|x|<n-1$,
 \begin{equation}
B_n(x)= {1\over \lambda} \sum_{i=1}^{d_x} { B_n(x_i) \over
1+ B_n(x_i)}, \quad \Gamma_n(x)= {1\over \lambda} \sum_{i=1}^{d_x} {
1+  \Gamma_n(x_i) \over 1+ B_n(x_i)}.\label{eq:11}
\end{equation}

Notice that we could define $\Gamma(x) : =
\lim_{n\to \infty} \frac{\Gamma_n(x)}n$, such that
\begin{equation}
  \label{eq:10}
  \Gamma(x):= {1\over \lambda} \sum_{i=1}^{d_x} {\Gamma(x_i) \over 1+ B(x_i)}
\end{equation}
As $\alpha\to 0$ we can show
that $\Gamma(x) \to aY(x)$ for some constant $a>0$. The problem is that
in the limit as $n\to\infty$ we loose information on the value of $a$
(that should be $\mathbb E(d_o(d_o -1))/[m(m-1)]$).

In order to determine this constant we have to make a step back and
iterate the equations \eqref{eq:11} and, noticing that $\Gamma_n(x)=0$ for  all
$\vert x \vert =n-1$, we get that   \begin{equation}\label{Gamma3}
  \Gamma_n(o)= \sum_{r=1}^{n-1} 
 {1\over \lambda^r} \sum_{|u|=r} { 1 \over 1+ B_n(u_1)} \cdot\cdot
 \cdot {1\over 1+ B_n(u_{r-1})}  {
1  \over 1+ B_n(u_r)} := \sum_{r=1}^{n-1} \Phi_n(r) , \end{equation}

\noindent where $\{u_0, ...., u_r\}$ is the shortest path relating
the root  $o$ to $u$ [$u_0=o, |u_1|=1, ..., |u_r|=r$]. Note
that
$B_n(u_1), ..., B_n(u_r)$ are correlated.

Observe that $\Phi_n(r) \le e^{\alpha r} W(o,r)$, consequently $\mathbb
E(\Phi_n(r)) \le e^{\alpha r}$. Since $\alpha>0$ it is hard to control
the limit of $\Gamma_n(o)$. The aim is to analyze  the asymptotic behavior of
$\e(\Phi_n(r))$ as $n \to \infty$ and $r\le n$, which will be done in
the following two subsections: in the next first subsection we will
give a useful representation of $\e(\Phi_n(r))$ based on a spine
random walk, whereas in  the second subsection  we make use of an
argument from renewal theory and establish the limit of
$\e(\Phi_n(r))$ when $r, n \to \infty$ in an appropriate way.

\subsection{Spine random walk representation of $\e(\Phi_n(r))$}

Let $\Omega$ denote the space of rooted trees with no leaves.
Denote by $\widetilde \Omega_T$  the space  of trees with a marked infinite ray
$\Ray=(u^*_n)_{n\ge0}$, with $u_0^*=o$ [$\widetilde \Omega_T$ is topologically equal to $ \Omega_T$]. 
Unlike the setup used in Section \ref{sec-2}, where 
e.g. $u_1^*$ was considered a parent of $o$,
we now
redefine the notion of descendant in $\widetilde 
\Omega_T$. Namely, for $x\in {\cal T}$,
$x\neq o$,
the parent of $x$, denoted ${\buildrel
\leftarrow \over x}$, is the unique vertex on the geodesic connecting
$x$ and $o$ with $|{\buildrel
\leftarrow \over x}|=|x|-1$. In this section,    for any $\vert v \vert  <  n$, we define  the normalized progeny of  $v$ at level $n$ as   $M_n(v):=|\{w: |w|=n,  w \mbox{ descendant of } v\}|/m^{n- \vert v\vert}$, and $M_n(v)=1$ if $\vert v \vert =n$.  We also write $M_n=M_n(o)$.


According to \cite{L97},
on the space $\widetilde \Omega_T$
we may construct a probability $\q$ such that  the marginal of 
$\q$ on the space of trees $\Omega$
satisfies $$ {d \q \over d \p} \vert_{ {\cal F}_n}:= M_n,  \qquad n\ge1, $$ 
where  ${\cal F}_n$ is the $\sigma$-field generated by the first
$n$ generations of $\omega$ and $\p$ denotes the Galton--Watson law.
Due to $p_0=0$ and our tail assumptions, we have that
$M_\infty>0$ a.s. and moreover,  
$$ \q\Big( u_n^*= u \big \vert {\cal F}_n \Big)=  { 
1\over M_n m^n}, \qquad \forall \, \vert u \vert =n.$$  
Under $\q$,  $ d_{u^*_n}$  has the size-biased distribution associated
with $\{p_k\}$, that is $\q(d_{u^*_n}=k)=kp_k/m$,
$u_{n+1}^*$ is uniformly chosen among the children of $u^*_n$ and,
for   $v \not= u_{n+1}^*$ with ${\buildrel
\leftarrow \over v} =u_n^*$, 
the sub-trees $\T(v)$ rooted at $v$, are i.i.d.  and have a Galton--Watson law.



For any $ 0\le j < n$, we define $a_j^{(n)}:={1\over \lambda}  \sum_{ v \not = u^*_{j+1}, 
{\buildrel \leftarrow \over v}=  u^*_j } \beta_n(v)$.  Note that under $\q$, the family $\{a_j^{(n)}\}_{0\le j < n}$ are independent and each $a_j^{(n)}$ is distributed as ${1\over \lambda} \sum_{k=1}^{ d^*-1} \beta_{n-j -1}^{(k)}$, where $d^*$   has the size-biased distribution associated with $\{p_k\}$, $(\beta^{(k)}_l, l\ge 1)$ are i.i.d copies of $(\beta_l (o), l\ge1)$ and independent of $d^*$.  We extend $a^{(n)}_j$  to  all $ j \in \z \cap (-\infty, n-1]$ by letting the family $\{a^{(n)}_j, n>j\}_{j \in \z}$ be independent  (under $\q$) and such that for each $j$, $\{a^{(n)}_j, n >j\}$ is distributed as 
$\{{1\over \lambda} \sum_{k=1}^{ d^*-1} \beta_{n-j -1}^{(k)}, n >j\}$.  We  naturally define   $a^{(\infty)}_j$ as the limit of $a^{(n)}_j$ as $n \to \infty$.  In particular for $j\ge 0$, $a_j^{(\infty)}:={1\over \lambda}  \sum_{ v \not = u^*_{j+1}, 
{\buildrel \leftarrow \over v}=  u^*_j } \beta(v),$ and  each  $a_j^{(\infty)}$   is  distributed as ${1\over \lambda} \sum_{k=1}^{ d^*-1} \beta ^{(k)}$ with $(\beta^{(k)})_{k\ge1}$ i.i.d copies of $\beta(o)$, independent of $d^*$.  


The main result of this subsection is the following representation for $\e(\Phi_n(r))$:

\begin{prop}\label{p:31}  We may define a random walk $(S_\cdot, P )$ on $\z$,  independent of the Galton-Watson tree $\omega$ and of the family $(a^{(n)}_j)_{j < n}$,  with  step distribution $ P (S_i-S_{i-1}= 1)={ \lambda \over \lambda+m ^2}$ and  $P (S_i-S_{i-1}=-1)={m^2 \over \lambda+m ^2}$,   $ \forall  i\ge1,$  such that for any $1\le r \le n$, 
 \begin{equation}\label{ephinr55bis} \e(\Phi_n(r)) = \q \left[  { Z_n(r) \over M_{ n-r} } \, \right] ,  \end{equation} where $$  Z_n(r):= E_{0, \omega} \Big( {\bf 1}_{( \tau_{S}(-r) < \tau_{S}(n-r))}  \prod_{i=0}^{\tau_S(-r)-1}  f_{n-r}(S_i)\Big)  ,  \qquad 1\le r \le n, $$ 
with $E  _{0, \omega}$  the expectation with respect to the random walk $S$ starting from  $0$,   and        \begin{equation} f_n(x):=  {m^2+\lambda   \over m(1+ \lambda + 
\lambda a^{(n)}_{x})}  , \qquad x< n .  \label{def-fn} \end{equation}  
 \end{prop}

Before entering into the proof of Proposition \ref{p:31}, we mention that the random walk $(S_\cdot, P )$ may find its root in the following lemma:

\begin{lemma}[Spine random walk] \label{L:zjbis}
 Let $n>  k\ge 2$. Let $b_{j+1}>0$ and $a_j \ge 0$ for all
 $0\le j < n$. Define $(z_j)_{0\le j\le n}$ by $z_n=0$ and
 $$
 z_j
 :=
 {1\over 1+  a_j + b_{j+1} (1-z_{j+1}) }, \qquad 0\le j\le n-1.
 $$
Let  $ (S_m)$ be  a Markov chain on $\{0, 1, ..., n\}$   
with probability transition $\widetilde P(S_m= j+1\vert S_{m-1}= j)= { b_{j+1}
\over 1+ b_{j+1}}$ and  $\widetilde P(S_m= j-1\vert S_{m-1}= j)= 
{ 1 \over 1+ b_{j+1}}$, and denote by $\widetilde P_r$ the law of the chain
$(S_m)$ with $S_0=r$.  
Then, for any $1\le r < n$,
$$    \prod_{j=1}^{r} z_j = \widetilde E_r \Big( {\bf 1}_{( \tau_{S}(0) < \tau_{S}(n))}
    \prod_{j=1}^{n-1}  \Big( {1+b_{j+1} \over 1+b_{j+1} +a_j}\Big)^{    L_{\tau_{ S}(0)}^ j }\Big),$$ with $\tau_S(x):=\inf\{ j\ge 1: S_j=x\}$ the first hitting time of $S$ at $x$ and $L_m^x:= \sum_{i=0}^{m-1} 1_{(S_i=x)}$ is the local time at $x$.
    \end{lemma}

 Lemma \ref{L:zjbis} can be proved  exactly as     in  \cite[Appendix]{FHS}, by using (A.3) and the construction of the random walk therein.  
We omit the details. Thanks to the spine random walk,  studying  $\e(\Phi_n(r))$ reduces  to a problem of one-dimensional random walk $S_\cdot$ in random medium (whose laws are determined by that of $\beta_n(\cdot)$); we solve the latter problem   by   using the regeneration times for the transient walk $S$ and the renewal theorem. 

\medskip

{\noindent \bf Proof of Proposition \ref{p:31}:}  Observe that $\beta_n(x)= { B_n(x) \over 1+ B_n(x)}$
and
$$ \Phi_n(r)= { 1\over \lambda^r } \sum_{ \vert u \vert =r} (1- \beta_n(u_1))
\cdots (1- \beta_n(u_r)) .$$

By the change of measure, we have for any $F\ge0$,
$$ \e \left[\sum_{ \vert u \vert =n} F( \beta_n(u_1), d_{u_1}, ...., \beta_n(u_n), d_{u_n}) \right]= m^n \q \left[ F( \beta_n(u^*_1), d_{u_1^*}, ..., \beta_n(u^*_n), d_{u^*_n})\right].$$

\noindent It follows that  for
$r<n$, \begin{eqnarray*} && m^{r} \,    \q \Big( (1- \beta_n(u^*_1))\cdots
(1- \beta_n(u^*_r)) {1\over M_n(u^*_r)} \Big) \\ 
    &=& m^{-(n-r)} \e\left[ \sum_{ \vert v \vert =n} (1- \beta_n(v_1))\cdots
(1- \beta_n(v_r)) { 1\over M_n(v_r)} \right] \\
    &=& \e \left[\sum_{ \vert u \vert =r} (1- \beta_n(u_1))\cdots(1- \beta_n(u_r))  \right], \end{eqnarray*}

\noindent where the term $m^{-(n-r)}{1\over M_n(v_r)}$ disappears when one takes the sum over $\vert v \vert =n$ by keeping $v_r=u$.   It follows  that  \begin{equation}\label{ephinrbis} \e(
\Phi_n(r) )= { m^r \over \lambda^r} \q \Big( (1- \beta_n(u^*_1))\cdots
(1- \beta_n(u^*_r))  {1\over M_{n }(u^*_r)} \Big), \end{equation}

\noindent and exactly the same as \eqref{Gamma3}, by iterating the equations on $B_n$, we get that
for any $r \le n-1$,
$$  B_n(o)= {1\over \lambda^r} \sum_{|u|=r} { 1 \over 1+ B_n(u_1)} \cdot\cdot
 \cdot {1\over 1+ B_n(u_{r-1})}  { B_n(u_r) \over 1+ B_n(u_r)}.$$
 Hence, 
\begin{equation}\label{ebn1bis} \e( B_n(o))= 
	{ m^r \over \lambda^r} \q \Big( (1- \beta_n(u^*_1))\cdots
	(1- \beta_n(u^*_{r-1}))  \beta_n(u^*_r)  {1\over M_{n }(u^*_r)} \Big)  . \end{equation}

Note that
$$ \beta_n(u^*_j)= {  \beta_n(u^*_{j+1})+ \sum_{ v \not = u^*_{j+1}, 
{\buildrel
\leftarrow \over v}=  u^*_j } \beta_n(v) \over
\lambda + \beta_n(u^*_{j+1})+ \sum_{ v \not = u^*_{j+1}, 
{\buildrel
\leftarrow \over v}=  u^*_j } \beta_n(v)} = {  \beta_n(u^*_{j+1})+  \lambda \, a_j^{(n)} \over
\lambda + \beta_n(u^*_{j+1})+    \lambda \, a_j^{(n)}} , \qquad  \forall j < n, $$

\noindent with  $\beta_n(u^*_n)=1$,  and $$ M_{n }(u^*_j)= { 1\over m} \sum_{ v \not = u^*_{j+1}, 
{\buildrel
\leftarrow \over v}=  u^*_j }  M_{n }(v) + {1\over m} M_{n }(u^*_{j+1}),   \qquad  \forall j < n, $$

\noindent with $M_n(u^*_n)=1$. Under $\q$, for such $|v|=j+1$, $(\beta_n(v), M_{n }(v)) $ are i.i.d. and distributed as $(\beta_{n-j-1}(o), M_{n-j-1}(o))$ (under $\p$). 


We can represent $1- \beta_n(u^*_j)$ as the probability for a one-dimensional  random walk in a random
environment (RWRE) with cemetery point, 
starting from $j$,  to hit $j-1$ before $n$.  In fact, 
applying Lemma \ref{L:zjbis} to $a_j=a^{(n)}_j $ and 
$b_{j+1} = { 1\over \lambda}$, we see that 
\begin{eqnarray*} \prod_{j=1}^r (1-  \beta_n(u^*_j)) &= & 
	\widetilde E_{r, \omega} \Big( {\bf 1}_{( \tau_{S}(0) < \tau_{S}(n))}
    \prod_{j=1}^{n-1}  \Big( {1+\lambda   \over 1+ \lambda + \lambda a^{(n)}_j}\Big)^{    L_{\tau_{ S}(0)}^ j }\Big)
    \\&=& \widetilde  E_{r, \omega} \Big( {\bf 1}_{( \tau_{S}(0) < \tau_{S}(n))}
    \prod_{i=0}^{\tau_S(0)-1}   {1+\lambda   \over 1+ \lambda + \lambda a^{(n)}_{S_i}} \Big) ,  \end{eqnarray*} where $(S_i)_{i\ge0}$     is a random walk on $\z$  with step distribution  
    $\widetilde  P(S_i- S_{i-1}=1)= { 1\over 1+\lambda}$ and 
    $\widetilde P(S_i- S_{i-1}=-1)= { \lambda\over 1+\lambda}$ for $i\ge1$, and
    the expectation $\widetilde E_{r, \omega}$ is taken with respect to $(S_m)$ with $S_0=r$.

Define  the  probability $P $ with $$ { d P  \over d \widetilde P} \Big\vert_{\sigma\{S_0, ..., S_n\}}=   \left( { \lambda\over m}\right)^{S_n- S_0 }\, \left( { m(1+\lambda) \over m^2 + \lambda}\right)^n ,   \qquad n\ge0. $$ 
Under $P$, the random walk $\{S_m\}$ has the properties stated 
in the statement of Proposition
\ref{p:31}.
Further, 
$$ {m^r \over \lambda^r} \prod_{j=1}^r (1-  \beta_n(u^*_j))  = E  _{r, \omega} \Big( {\bf 1}_{( \tau_{S}(0) < \tau_{S}(n))}
    \prod_{i=0}^{\tau_S(0)-1}  f_n(S_i)\Big)  := \widetilde Z_n(r) .$$

\noindent With the notation of $\widetilde Z_n(r)$,   we get that \begin{equation}\label{ephinr5bis}  \e \left(\Phi_n(r)\right)=   \q \left[ {\widetilde  Z_n(r)    \over M_n(u^*_r)} \right]. \end{equation}

Observe that for any   $r< n$, under $\q$, $(f_n(x+r), M_n(u^*_r))_{x \le n-r }$ has the law as  $(f_{n-r}(x ), M_{n-r})_{x \le n-r }$. This invariance by linear shift  and   (\ref{ephinr5bis})  yield Proposition \ref{p:31}. $\Box$

    \medskip 
    
    We end this subsection by the following remark:
    
\begin{remark} \label{r:33} With the the same notations as in Proposition \ref{p:31}, we have \begin{eqnarray}  \e (B_n(o)) & \le &  { m\over \lambda} \q \left[ { Z_n(r-1)\over M_{n-(r-1)}(u^*_1) }\right]  , \label{ebn4bis}
    \\ \e (B_n(o)) &\ge & { m\over \lambda} \q \left[ { Z_n(r-1)\over M_{n-(r-1)}(u^*_1)} \, { a_1^{(n-r+1)} \over 1+ a_1^{(n-r+1)}}   \right]  . \label{ebn5bis}\end{eqnarray}
\end{remark}
    
    {\noindent \bf Proof of Remark \ref{r:33}:}  
In the same way which leads to \eqref{ephinr5bis}, we get from (\ref{ebn1bis})   that \begin{eqnarray}  \e(B_n(o)) &= &{ m ^r \over \lambda^r} \, \q \Big( ( 1- \beta_n(u^*_1)) \cdot \cdot \cdot (1- \beta_n(u^*_{r-1})) \beta_n(u^*_{r})  \, {1\over M_{n }(u^*_r)} \Big) \nonumber\\
	&=& { m ^r \over \lambda^r} \, \q \Big( \Big[\prod_{i=1}^{r-1} (1- \beta_n(u^*_{i})) - \prod_{i=1}^{r } (1-\beta_n(u^*_{i})) \Big]   \, {1\over M^{*, n}_r } \Big)  \nonumber
    \\&=& { m  \over \lambda } \, \q \Big(   {\widetilde Z_n(r-1)        \over M_{n }(u^*_r) } \Big) - \q \Big(   { \widetilde Z_n(r )       \over M_{n }(u^*_r) } \Big),  \label{ebn2bis}\end{eqnarray}

\noindent  giving the upper bound in \eqref{ebn4bis} after a linear shift  at $r-1$ for the above term with $m\over \lambda$. On the other hand,  $$  \beta_n(u^*_r)= { \beta_n(u^*_{r+1}) +\lambda a^{(n)}_r   \over \lambda +\beta_n(u^*_{r+1}) +\lambda a^{(n)}_r  }  \ge  {   a^{(n)}_r  \over 1+ a^{(n)}_r }, $$

 \noindent hence \begin{eqnarray}  \e(B_n(o)) &\ge  &{ m ^r \over \lambda^r} \, \q \Big(( 1- \beta_n(u^*_1)) \cdot \cdot \cdot (1- \beta_n(u^*_{r-1})) {   a^{(n)}_r  \over 1+ a^{(n)}_r }  \,  {1\over M_{n }(u^*_r)} \Big) \nonumber\\&=& { m   \over \lambda } \, \q \Big( \widetilde  Z_n(r-1)     {   a^{(n)}_r  \over 1+ a^{(n)}_r }  \,  {1\over M_{n }(u^*_r)} \Big), \label{ebn3bis}\end{eqnarray}

\noindent yielding the assertions in Remark \ref{r:33} after the shit at $r-1$. $\Box$

\subsection{An  argument based on renewal theory}

The main result is Lemma \ref{lem-Udist} which evaluates  the limit of $\e(\Phi_n(r))$ and in turn gives the velocity representation in Theorem \ref{theo-3.4}.    The analysis is based on the
use of a 
renewal structure in the representation of 
Proposition \ref{p:31}.  Under $P $, $(S_i)$ drifts to $-\infty$.  
Denote by $(R_0:=0)<
R_1< R_2 <...$ the regeneration times for $(S_i)$, that is 
$R_i=\min\{n>R_{i-1}: \{S_i\}_{j=0}^n\cap \{S_j\}_{j>n}=\emptyset\}$.   
The sequence $\{S_{j+R_i}- S_{R_i}, 0\le j \le R_{i+1}-R_i\}_{i\geq 1}$ is clearly i.i.d and has as
common distribution that of $\{S_j,  0\le j \le R_1\}$ conditioned on $\{\tau_S(1)=\infty\}$.  Further, 
because 
$$E  (S_{i+1}-S_i)=\frac{\lambda-m^2}{\lambda+m^2}\leq -\frac{m-1}{m+1}\,,
$$
it is straightforward to check that there exists a constant $\kappa>0$, 
independent of $\alpha$, so that
\begin{equation}
	\label{reg-exp}
	E  (e^{\kappa R_1})<\infty\,,\quad
	E  (e^{\kappa(R_2-R_1)})<\infty\,.
\end{equation}

Define  
$$\zeta_j:= \prod_{i=R_{j-1}}^{R_j-1} f_\infty(S_i), \qquad j\ge 1, $$ where \begin{equation}
	\label{eq-fdef}
	f_\infty(x):=  {m^2+\lambda   \over m(1+ \lambda + \lambda a^{(\infty)}_{x})} , \qquad x \in \z\,, 
\end{equation}

\noindent and  for $x \in \z$,  $a_x^{(\infty)}$   are i.i.d. and  are distributed as ${1\over \lambda} \sum_{k=1}^{ d^*-1} \beta ^{(k)}$ with $(\beta^{(k)})_{k\ge1}$ i.i.d copies of $\beta(o)$, independent of $d^*$.

An important observation is  that under $\q \otimes P $, $(\zeta_j, j\ge 2)$ are i.i.d. and independent of $\zeta_1$.  Define further 
\begin{equation}
	\label{eq-hdef}
	h(y):=   \q \otimes P     \left[  \prod_{i=0}^{\tau_S(-y) -1}f_\infty(S_i)  1_{( \tau_S(-y) \le R_1 )} \, \Big \vert \, \tau_S(1)= \infty\right], \qquad y \ge 1. 
\end{equation}

\noindent We extend the definition of $h$ to $\z$ by letting $h(y):=0$ if $y \le 0$.

\begin{lemma} 
Assume that
  \begin{eqnarray} && \sum_{y=1}^\infty h(y) < \infty, \label{L:hyp1} \\ 
	  &&\q \otimes P     
	  \left[{ \zeta_1\over M_\infty(u^*_1)}+ \zeta_1+  | S_{R_2} - S_{R_1}| \, 
	  \zeta_2\right] <\infty \label{L:hyp2}\,. 
\end{eqnarray}
\noindent Then,
$$ \q \otimes P \left[ \zeta_2\right]=1.$$
\end{lemma}
We will see below, see Lemma \ref{lem-tildeqreg}, that \eqref{L:hyp1} and
\eqref{L:hyp2} both hold for all $\alpha>0$ small enough.

\bigskip
{\noindent\bf Proof:} Almost surely, $\beta_n(x) \downarrow \beta(x)$. Then for a fixed $r$, almost surely,  $$   Z_n(r) \to  Z_\infty(r):=  E  _{0, \omega} \Big(
    \prod_{i=0}^{\tau_S(-r)-1}  f_{\infty}(S_i)\Big) , $$ and $  M_n(u^*_1)  \to M_\infty(u^*_1)$. Applying Fatou's lemma in the expectation under $\q$ in (\ref{ebn5bis}), we get that for any $r$,  \begin{eqnarray*} {\lambda\over m} \e(B) &\ge &   \q \left[ { Z_\infty(r-1)\over M_\infty(u^*_1) } \,  {   a_1^{(\infty)} \over 1+ a_1^{(\infty)}  }  \right]
     \\&=&     \q \otimes P \left[ \prod_{i=0}^{\tau_S(1-r)-1}  f_{\infty}(S_i)   \, {1  \over M_\infty(u^*_1) } \,  {   a_1^{(\infty)} \over 1+ a_1^{(\infty)}  }   \right] . \end{eqnarray*}

We can not directly let $r \to \infty$ inside the above expectation, so we decompose this expectation by the regeneration times $0<R_1<R_2<...$. Write $$\zeta'_1:=     {\zeta_1  \over M_\infty(u^*_1) } \,  {   a_1^{(\infty)} \over 1+ a_1^{(\infty)}  }   .$$

\noindent Then    
 \begin{eqnarray*}  {\lambda\over m} \e(B)  & \ge  &\sum_{k=2}^\infty   \q \otimes P   \left[   1_{( R_k < \tau_S(1-r)\le R_{k+1})} \, \zeta'_1\, \prod_{i=R_1}^{R_k-1} f_\infty(S_i) \,  \prod_{i=R_k}^{\tau_S(1-r)-1} f_\infty(S_i)\right]
    	\\ &=&  \sum_{k=2}^\infty   \q \otimes P   \left[   1_{( R_k < \tau_S(1-r))} \, \zeta'_1\, \prod_{i=R_1}^{R_k-1} f_\infty(S_i) \,  h( r-1+ S_{R_k}) \right]  ,  \end{eqnarray*}

\noindent by using the Markov property of $S$  at $R_k$.   Observe that  $(\zeta_j, S_{R_j}- S_{R_{j-1}})_{j\ge 2} $ are i.i.d. under the annealed measure $\q \otimes P $, and are independent of $(\zeta'_1, S_{R_1})$. By replacing $r-1$ by $r$, we get that for any $r$,  \begin{equation}\label{eb6bis}\sum_{k=2}^\infty   \q \otimes P   \left[   1_{( S_{R_k} > -r)} \, \zeta'_1\, \prod_{j=2}^k  \zeta_j \,  h( r + S_{R_k}) \right] \le  {\lambda\over m} \e(B).\end{equation}

Now, we claim that \begin{equation}\label{<=1bis} \q \otimes P  \left[ \zeta_2\right] \le 1. \end{equation}

	To prove \eqref{<=1bis}, we
assume
that $a:=\q \otimes P  \left[ \zeta_2\right] >1$ and
show that
it  leads to a contradiction with (\ref{eb6bis}). Toward this end,
define a distribution $U$ on $ \z_+$ by  $$ U(x):= {\q \otimes P  \left[ 1_{(S_{R_2}- S_{R_1}= -x)} \zeta_2\right] \over \q \otimes P  \left[ \zeta_2\right]}, \qquad x \ge 0.$$

Then (\ref{eb6bis}) becomes 
\begin{eqnarray*} {\lambda\over m} \e(B) &\ge & \sum_{k=2}^\infty  \, a^{k-1} \,   \q \otimes P   \left[   1_{( S_{R_1} > -r)} \, \zeta'_1\, \sum_{x=0}^{ r+S_{R_1}}    h( r + S_{R_1}-x ) U^{\otimes (k-1)} (x) \right]
    \\&\ge& a^{l-1}\, \sum_{k=l}^\infty  \,     \q \otimes P   \left[   1_{( S_{R_1} > -r)} \, \zeta'_1\, \sum_{x=0}^{ r+S_{R_1}}    h( r + S_{R_1}-x ) U^{\otimes (k-1)} (x) \right], \end{eqnarray*} for any $l\ge 2$. 
    
    Since $\sum_{k=1}^{l-1}  \,     \q \otimes P   \left[   1_{( S_{R_1} > -r)} \, \zeta'_1\, \sum_{x=0}^{ r+S_{R_1}}    h( r + S_{R_1}-x ) U^{\otimes (k-1)} (x) \right] \to 0$ as $r \to \infty$ [by the dominated convergence under (\ref{L:hyp1}) and the integrability of $\zeta_1' \le { \zeta_1\over M_\infty(u^*_1)}$ under (\ref{L:hyp2})], we get that
    for 
    any fixed $\ell$,  \begin{eqnarray*}  a^{1-l}\, {\lambda\over m} \e(B) &\ge & \sum_{k=1}^\infty  \,     \q \otimes P   \left[   1_{( S_{R_1} > -r)} \, \zeta'_1\, \sum_{x=0}^{ r+S_{R_1}}    h( r + S_{R_1}-x ) U^{\otimes (k-1)} (x) \right] + o(1) \\
    &=& \q \otimes P   \left[    \zeta'_1  \right] {\sum_{x=0}^\infty h(x) \over \sum_{x=0}^\infty x U(x)} +o(1), \qquad r \to \infty,\end{eqnarray*}

    \noindent by applying the renewal theorem \cite[pg. 362]{Fel},
    using \eqref{L:hyp2}.
    Thus   we get some constant $C>0$ such that ${\lambda\over m} \e(B) \ge a^{l-1} C$ for any $\ell \ge2$, which is  impossible if  $a>1$. Hence we proved (\ref{<=1bis}).

It remains to show \begin{equation}\label{>=1bis} \q \otimes P  \left[ \zeta_2\right] \ge 1. \end{equation}

The proof of this part is similar, 
we shall use (\ref{ebn4bis}) instead of (\ref{ebn5bis}).
Set
$$\bar f_{S}(r):=
\prod_{i=0}^{\tau_S(-r)-1}  f_{\infty}(S_i) 
\,.$$
Since $f_\infty(x) \ge f_\ell(x)$ for any $\ell$, we get that $$ {\lambda\over m}  \e(B_n(o)) \le \q \otimes P   \left[    {\bf 1}_{( \tau_{S}(1-r) < \tau_{S}(n-(r-1)))}
\bar f_{S}(r-1)
    \, { 1\over  M_{n-r+1}(u^*_1) }\right] .$$

\noindent 
Taking $r=n$ gives that  $$ {\lambda\over m} 
\e(B_n(o)) \le \q \otimes P   \left[    {\bf 1}_{( \tau_{S}(1-n) < \tau_{S}(1))}\bar f_S(n-1)\right]\,.$$

\noindent 
Since $\e(B) \le \e(B_n(o))$, we obtain that for any $n$,  
\begin{eqnarray*} {\lambda\over m}  \e(B ) &\le& 
	\q \otimes P   \left[ 
	{\bf 1}_{( \tau_{S}( -n) < \tau_{S}(1))}\bar f_S(n)\right]\\
    &\le & \q \otimes P   \left[    
    {\bf 1}_{( \tau_{S}(-n) \le  R_1,\,  \tau_{S}( -n) < \tau_{S}(1))  )} \bar f_S(n) \right] +
    \sum_{k=1}^\infty \q \otimes P   \left[   
    {\bf 1}_{( R_k < \tau_{S}(-n) \le  R_{k+1})} \bar f_S(n)\right].
\end{eqnarray*}

\noindent
By the Markov property at $\tau_S(-n)$, the first term equals 
$$ {\q \otimes P   \left[    {\bf 1}_{( \tau_{S}(-n) \le  R_1)} 
\bar f_S(n)  {\bf 1}_{(      \tau_{S}(1) =\infty )}   \right]  \over P _{-n}(\tau_S(1)=\infty)} = 
{P  (\tau_S(1)=\infty)  \over P _{-n}(\tau_S(1)=\infty)} h(n) \to 0\,,$$
since $P _{-n}(\tau_S(1)=\infty) \ge c$ for some constant $c>0$ and $h(n) \to 0$.  
Then, recalling that $h$ vanishes at $\z_-$, we get
$$ {\lambda\over m} \e(B ) \le o(1)+ 
\sum_{k=1}^\infty \q \otimes P   
\left[    \zeta _1 \, \prod_{j=2}^k \zeta_j \, h ( n+ S_{R_k} )\right]\,, $$

\noindent with $\zeta_j$ and $h$ defined as before.  If $a:=\q \otimes P   \left[ \zeta_2\right] < 1$, then with the distribution $U(\cdot)$ introduced before, \begin{eqnarray*} \sum_{k=1}^\infty \q \otimes P   \left[    \zeta _1 \, \prod_{j=2}^k \zeta_j \, h ( n+ S_{R_k} )\right] &= &\sum_{k=1}^\infty  \, a^{k-1} \,   \q \otimes P   \left[     \, \zeta_1\, \sum_{x=0}^{ n+S_{R_1}}    h( n + S_{R_1}-x ) U^{\otimes (k-1)} (x) \right]
    \\&:=&  \sum_{k=1}^\infty  \, a^{k-1} \, b_k^{(n)}.\end{eqnarray*}

       \noindent Note
       that $\max_n b_k^{(n)} \le \q \otimes \widetilde 
       P  \left[     \, \zeta_1\right]\, 
       \sum_{x=0}^\infty h(x)\, U^{\otimes (k-1)} (x)$, and that, due
       to \eqref{L:hyp1},
       $ \lim_{n \to\infty} b_k^{(n)} = 0$.
      The dominated convergence theorem
      then 
      implies that  $\sum_{k=1}^\infty  \, a^{k-1} \, b_k^{(n)} \to 
      0$ which in turn yields ${\lambda\over m}  \e(B ) \le o(1)$,
       a contradiction. Thus  $\q \otimes P  
       \left[ \zeta_2\right] \ge 1$. This
       completes the proof of the lemma.
$\Box$

\bigskip

\begin{lemma} \label{lem-Udist}
	Assume \eqref{L:hyp1}, 
	\eqref{L:hyp2} and that for some $p>1$, \begin{equation} \label{hyp-zeta2}  \q \otimes P  ( (\zeta_2 )^p)< \infty. \end{equation}  Furthermore, we assume  that 
		\begin{equation} \label{unif}  \mbox{ under $ \q \otimes P$,  the family $\{ {\zeta_1 \over M_k}\}_{k\ge1}$ is uniformly integrable}, \end{equation} and 	that \begin{equation}  	\label{L:hyp3} \lim_{  r\to\infty}\, \sup_{n\ge r}  \q \otimes P   \left[
    {\bf 1}_{( \tau_{S}(-r) \le R_1)} \prod_{i=0}^{\tau_S(-r)-1}
f_{\infty}(S_i)   \, {1\over M_{n-r}} \, \right] =0, \end{equation}  where as before,  $R_1$ is the first regeneration time for $S$ under $P $. Then,
	for any $\varepsilon>0$,
\begin{equation} \label{const} \lim_{n\to\infty} \max_{\varepsilon n \le r \le (1-\varepsilon) n} \left| \e(\Phi_n(r))-  {  \q \otimes P   \left({ \zeta_1\over M_\infty} \right) \, \sum_{y\ge1} h(y) \over  \q \otimes P   \left( \zeta_2 \vert S_{R_2} - S_{R_1}\vert \right) } \right|=0. \end{equation} Moreover, \begin{equation}\label{ephi99}  \sup_{r \ge 1, \, n \ge r  } \e ( \Phi_n(r)) < \infty. \end{equation}
  \end{lemma}
{\noindent \bf Proof:}    We  split the proof of \eqref{const} into the following upper and lower bounds:   \begin{eqnarray}   \label{Udist:upper} \limsup_{r \to\infty, \, n-r \to \infty}  \e(\Phi_n(r) ) &\le& {  \q  \otimes P   \left( { \zeta_1\over M_\infty}\right) \, \sum_{y\ge 1} h(y) \over  \q \otimes P   \left( \zeta_2 \vert S_{R_2} - S_{R_1}\vert \right) }   , \\
  \liminf_{n\to\infty} \min_{ \varepsilon n \le r \le (1-\varepsilon) n}  \e (\Phi_n(r)) &\ge& {  \q  \otimes P   \left( { \zeta_1\over M_\infty}\right) \, \sum_{y\ge 1} h(y) \over  \q \otimes P   \left( \zeta_2 \vert S_{R_2} - S_{R_1}\vert \right) } . \label{Udist:lower} \end{eqnarray}

\medskip
{\noindent \it Proofs of (\ref{ephi99}) and (\ref{Udist:upper}) :}   Let us  introduce the notation: for $\ell \ge 1$,
$$ \zeta_j(\ell):= \prod_{i=R_{j-1}}^{R_j-1} f_\ell(S_i), \qquad j\ge 1.$$
By (\ref{ephinr55bis}), $$  \e  (\Phi_n(r)) = \q \otimes P   \left[
    {\bf 1}_{( \tau_{S}(-r) < \tau_{S}(n-r))} \prod_{i=0}^{\tau_S(-r)-1}
f_{n-r}(S_i)   \, {1\over M_{n-r} } \, \right] .$$

Noticing that  $  {\bf 1}_{( \tau_{S}(-r) < \tau_{S}(n-r))}= 1_{(\tau_S(-r) \le R_1 \wedge \tau_S(n-r))}+ \sum_{k=1}^\infty 1_{(R_1 < \tau_S(n-r), \, R_k < \tau_S(-r) \le R_{k+1})}$, we get   \begin{eqnarray}   \e  (\Phi_n(r))&=& I_{(\ref{ephinr7bis})} (0)+\sum_{k=1}^\infty  I_{(\ref{ephinr7bis})} (k) ,  \label{ephinr7bis}\end{eqnarray}

\noindent where   
\begin{eqnarray*}  I_{(\ref{ephinr7bis})} (0)  &:= &     \q \otimes P   \left[
    {\bf 1}_{( \tau_{S}(-r) < R_1 \wedge \tau_{S}(n-r))} \prod_{i=0}^{\tau_S(-r)-1}
f_{n-r}(S_i)   \, {1\over M_{n-r}  } \, \right]  , \\
I_{(\ref{ephinr7bis})} (k) &:=&    \q \otimes P   \left[   1_{(R_1 < \tau_S(n-r), \,  R_k < \tau_S(-r)\le R_{k+1})}
   { \zeta_1(n-r) \over  M_{n-r} }\, \prod_{j=2}^{k} \zeta_j(n-r) \,  \prod_{i=R_k}^{\tau_S(-r)-1} f_{n-r}(S_i)\right]   \nonumber
    	\\&=&    \q \otimes P   \left[   1_{( R_1 < \tau_S(n-r), \, R_k < \tau_S(-r) )}
    { \zeta_1(n-r) \over M_{n-r} }\, \prod_{j=2}^{k} \zeta_j(n-r) \,   h_{n-r} ( r + S_{R_k}) \right]   ,   \end{eqnarray*}

\noindent  by the Markov property at $R_k$ (with convention $\prod_\emptyset \equiv 1$) and with  $$ h_{\ell}(y):=   \q \otimes P     \left[  \prod_{i=0}^{\tau_S(-y) -1}f_\ell(S_i)  1_{( \tau_S(-y) \le R_1 )} \, \Big \vert \, \tau_S(1)= \infty\right], \qquad y \ge 1, \, \ell\ge1. $$

 \noindent 
We also define $h_L(y):=0$ for all $y \le 0$. 
Since $f_\ell(x) \le f_\infty(x)$, $h_\ell(x) \le h(x)$, $\zeta_j(n-r) \le \zeta_j$ for any $j \ge 1$,  we have  $$ I_{(\ref{ephinr7bis})} (0) \le  \q \otimes P   \left[    {\bf 1}_{( \tau_{S}(-r) < R_1  )} \prod_{i=0}^{\tau_S(-r)-1}
f_{\infty}(S_i)   \, {1\over  M_{n-r} } \, \right] \, \to \, 0, \qquad \mbox{as } n>r \to \infty,  $$ 
by  \eqref{L:hyp3}.  Moreover, $I_{(\ref{ephinr7bis})} (0) $ is uniformly bounded over all $n\ge r \ge 1$, again by \eqref{L:hyp3}. Further,
\begin{eqnarray}  \sum_{k=1}^\infty   I_{(\ref{ephinr7bis})} (k)  &\le& \sum_{k=1}^\infty \q \otimes P   \left[   1_{( S_{R_1}>-r  )}
   { \zeta_1  \over M_{n-r}  }\, \prod_{j=2}^{k} \zeta_j  \,   h  ( r + S_{R_k}) \right]   \nonumber \\
   &=&   \sum_{k=1}^\infty    \q \otimes P     \left[    1_{( S_{R_1}>-r  )}
   { \zeta_1  \over  M_{n-r} } \, \sum_{x=0}^{ r + S_{R_1}} h( r+S_{R_1}- x) \, U^{\otimes (k-1)}(x) \right]  \nonumber \\
   &=&        \q \otimes P     \left[    1_{( S_{R_1}>-r  )}
   { \zeta_1  \over  M_{n-r} } \, \sum_{x=0}^{ r + S_{R_1}} h( r+S_{R_1}- x) \, u (x) \right] , \label{renewal2} \end{eqnarray}

\noindent where $u(x):= \sum_{k=1}^\infty U^{\otimes (k-1)}(x), x\ge 0, $ and $$ U(x):=   \q  \otimes P     \left[   1_{( S_{R_2} - S_{R_1} =  -x   )} \,  \zeta_2 \,    \right], \qquad x \ge 0,  $$ is a distribution by Lemma \ref{lem-Udist}.  By the renewal theorem, $$ u(x)   \to { 1 \over \sum y U(y)}= { 1 \over \q \otimes P   \left( \zeta_2 \vert S_{R_2} - S_{R_1}\vert \right) } := u(\infty)  , \qquad x \to \infty.  $$

\noindent  Recall \eqref{L:hyp1} that $\sum h(y) < \infty$. It follows    that as $r \to \infty$ and $n-r \to\infty$,  the  term $\left[ \cdot\cdot\cdot \right]$  in   (\ref{renewal2}) converges  almost surely to ${\zeta_1 \over M_\infty } \sum_{y=1}^\infty  h(y) u(\infty)$.  This  in view of  the uniform integrability \eqref{unif},  yield that  (\ref{renewal2}) converges  to  $$ {  \q  \otimes P   \left( { \zeta_1\over M_\infty}\right) \, \sum_{y\ge 1} h(y) \over  \q \otimes P   \left( \zeta_2 \vert S_{R_2} - S_{R_1}\vert \right) } , \qquad \mbox{ as } r \to\infty \mbox { and }  \, n-r \to \infty . $$

\noindent The estimate (\ref{Udist:upper})  follows.   Finally,  note that (\ref{renewal2})  is bounded by $$ \max_{x\ge 0} u(x) \, \sum_{y\ge1} h(y)\,   \q \otimes P     \left[     
   { \zeta_1  \over  M_{n-r} }  \right] \le c ,$$ for some constant $c>0$, uniformly over $n \ge r \ge 1$, again by \eqref{unif}.   Hence $\e( \Phi_n(r)) \le  I_{(\ref{ephinr7bis})} (0) + c $,   implying   \eqref{ephi99}.

\medskip
{\noindent\it Proof of (\ref{Udist:lower}) :}   The idea is to replace $\zeta_2(n-r)$ by $\zeta_2\equiv \zeta_2(\infty)$ in (\ref{ephinr7bis}).   Let $\ell = n-r$ and recall that $ f_\ell(x):=  {m^2+\lambda   \over m(1+ \lambda + \lambda a^{(\ell)}_{x})} ={m^2+\lambda   \over m(1+ \lambda + \sum_{k=1}^{d^*_x} \beta^{(x, k)}_{\ell-x-1} ) },  $  for $x < \ell$.
  Then $$ 0\le f_\infty(x)- f_\ell(x)=  {(m^2+\lambda ) \lambda ( a_x^{(\ell)}- a_x^{(\infty)})   \over m(1+ \lambda + \lambda a^{(\ell)}_{x}) \, (1+ \lambda + \lambda a^{(\infty)}_{x})}   = f_\infty(x) \, { \lambda ( a_x^{(\ell)}- a_x^{(\infty)})   \over   1+ \lambda + \lambda a^{(\ell)}_{x}  }  .$$

\noindent It follows that for any $j \ge2$, $$ \zeta_j(\ell) = \zeta_j\, \prod_{i= R_{j-1}}^{R_j -1}\, \left[ 1- { \lambda ( a_{S_i}^{(\ell)}- a_{S_i}^{(\infty)})   \over   1+ \lambda + \lambda a^{(\ell)}_{S_i}  }  \right]:= \zeta_j\, \times\, \Lambda_j(\ell).$$

Fix a large integer $L$.   Using (\ref{ephinr7bis}) and the fact that $h_\cdot(y)$ is nondecreasing for any $y$, we deduce that for all $n-r \ge L$ and any large constant $C>0$ [the constant $C$ will be chosen later on],  \begin{eqnarray*} \e(\Phi_n(r)) & \ge&   \sum_{k=1}^{C n}    \q \otimes P   \left[   1_{(R_1 < \tau_S(L), \,  R_k < \tau_S(-r) )}     {\zeta_1(L)\over  M_{n-r}  }\, \prod_{j=2}^{k} \zeta_j \,   \prod_{j=2}^{k} \Lambda_j(n-r) h_{L} ( r + S_{R_k}) \right].
\end{eqnarray*}

\noindent  The first step is to replace $\Lambda_j(n-r)$ by $1$,  then we have to check that   the error term is uniformly small:    \begin{equation}\label{err1bis}   I_{(\ref{err1bis})} :=  \sum_{k=1}^{C n}    \q \otimes P   \left[   1_{( R_1 < \tau_S(L), \, R_k < \tau_S(-r) )}     {\zeta_1(L)\over  M_{n-r} }\, \prod_{j=2}^{k} \zeta_j \,   \Big[1- \prod_{j=2}^{k} \Lambda_j(n-r) \Big] h_{L} ( r + S_{R_k}) \right] \to 0, 
\end{equation}

\noindent as $r \to \infty$ and $ \varepsilon n \le r \le (1- \varepsilon) n$. Let us postpone for the moment the proof of (\ref{err1bis}).   Going back to $\e(\Phi_n(r))$, we
obtain that
 \begin{eqnarray}&& \e(\Phi_n(r))  \nonumber \\& \ge&   \sum_{k=1}^{C n}    \q \otimes P   \left[   1_{( R_1 < \tau_S(L), \, R_k < \tau_S(-r) )}    { \zeta_1(L) \over M_{n-r}  } \, \prod_{j=2}^{k} \zeta_j \,    h_{L} ( r + S_{R_k}) \right] -  I_{(\ref{err1bis})}  \nonumber \\
 	&\ge&  \sum_{k=1}^{\infty}    \q \otimes P   \left[   1_{( R_1 < \tau_S(L), \,  R_k < \tau_S(-r) )}   { \zeta_1(L) \over M_{n-r}  } \, \prod_{j=2}^{k} \zeta_j \,    h_{L} ( r + S_{R_k}) \right] -  I_{(\ref{err1bis})} -   I_{(\ref{err1ter})} , \label{err1ter}
\end{eqnarray}

\noindent  with $$ I_{(\ref{err1ter})} := \sum_{k=Cn +1}^{\infty }    \q \otimes P   \left[   1_{( R_1 < \tau_S(L), \, R_k < \tau_S(-r) )}    { \zeta_1(L) \over M_{n-r} } \, \prod_{j=2}^{k} \zeta_j \,    h_{L} ( r + S_{R_k}) \right] .$$

If we can prove that  for a well-chosen $C$, $ I_{(\ref{err1ter})}$ goes to zero uniformly  for $r \to\infty$ and $\varepsilon n \le r \le (1-\varepsilon) n$, then  by applying  the renewal theorem ($L$ fixed, $r \to\infty$ and $n-r \to \infty$) to the   sum  in  (\ref{err1ter}), under  the uniform integrabiltiy \eqref{unif},  we get that $$ \liminf_{n\to\infty} \min_{n-r \ge \varepsilon n} \, \e(\Phi_n(r)) \ge   {  \q  \otimes P   \left( { \zeta_1(L) \over M_\infty}1_{( R_1 < \tau_S(L)  )} \right) \, \sum_{y\ge 1} h_L(y) \over  \q \otimes P   \left( \zeta_2 \vert S_{R_2} - S_{R_1}\vert \right) } . $$

\noindent Letting $L\to\infty$ gives the lower bound (\ref{Udist:lower}).

It remains to show that $ I_{(\ref{err1ter})} $ and  $ I_{(\ref{err1bis})} $ go to zero uniformly  for $r \to\infty$ and $\varepsilon n \le r \le (1-\varepsilon) n$.  We first deal with  $ I_{(\ref{err1ter})} $.  Let $h^*:=\max_{x\ge 0} h(x)$.  We have \begin{eqnarray*}   I_{(\ref{err1ter})}  &\le &  h^* \, \sum_{k=Cn +1}^{\infty }    \q \otimes P   \left[   1_{(   R_k < \tau_S(-r) )}    { \zeta_1(L) \over M_{n-r} } \, \prod_{j=2}^{k} \zeta_j  \right]   \\ 
	&\le &    h^* \, \sum_{k=Cn +1}^{\infty }    \q \otimes P   \left[       { \zeta_1(L) \over M_{n-r} } \,1_{(   S_{R_k} - S_{R_1} > -r  )} \prod_{j=2}^{k} \zeta_j  \right]   \\ 
	&=  &  h^* \,  \sum_{k=Cn +1}^{\infty }    \q \otimes P   \left[       { \zeta_1(L) \over M_{n-r} } \right]\,   \q \otimes P   \left[   \,1_{(   S_{R_k} - S_{R_1} > -r  )} \prod_{j=2}^{k} \zeta_j  \right] ,
\end{eqnarray*}

\noindent by the independence between $(  \zeta_1(L), M_{n-r} )$ and $(S_{R_k} - S_{R_1} , \zeta_j , j\ge 2)$.  Recall
that $ \q \otimes P   (\zeta_2)=1$.   Let  $\widehat P$  be a new probability measure defined by 
${d \widehat P \over d   \q \otimes P } = \zeta_2$, then under $\widehat P$, $S_{R_1} - S_{ R_k} $ is the sum of $k-1$ positive i.i.d. variables with mean $  \q \otimes P  ( \zeta_2 (S_{R_1} - S_{ R_2})) :=   a  \in (0, \infty)  $ by (\ref{L:hyp2}).   Taking $C:= {2 \over a}$.  Then by Cramer's bound, there exists some $c_0 >0$ such that  $$ \q \otimes P   \left[   \,1_{(   S_{R_k} - S_{R_1} > -r  )} \prod_{j=2}^{k} \zeta_j  \right] = \widehat P \Big( S_{R_1} - S_{R_k} < r \Big) \le e^{ - c_0 k}, $$ for any $k > C n$ and $r \le n$.  It follows $$  I_{(\ref{err1ter})}  \le h^* \,  \q \otimes P   \left[       { \zeta_1(L) \over M_{n-r} } \right]\,   \sum_{k=Cn +1}^{\infty }  e^{ - c_0 k} \to 0,$$   uniformly as $r \le n $ and $ r \to \infty$, since $\q \otimes P   \left[       { \zeta_1(L) \over M_{n-r} } \right] \le \q \otimes P   \left[       { \zeta_1  \over M_{n-r} } \right] \le C'$, with some constant $C'>0$,  thanks to the uniform integrability \eqref{unif}.

It remains to  check (\ref{err1bis}) [with $C:={2\over a}$ chosen before]. We first observe that $h_l(x) \le h(x)\le h^*$ and that \begin{eqnarray*} I_{(\ref{err1bis})}  &\le&  h^* \, \sum_{k=1}^{C n}    \q \otimes P   \left[      {\zeta_1 \over M_{n-r} } \, \prod_{j=2}^{k} \zeta_j \,   \Big(1- \prod_{j=2}^{k} \Lambda_j(n-r) \Big)  \right]  \\
	&=&  h^* \,  \sum_{k=1}^{C n}    \widehat E  \left[       {\zeta_1 \over M_{n-r} } \,   \Big( 1- \prod_{j=2}^{k} \Lambda_j(n-r) \Big)  \right] \\
	&=& h^* \,   \sum_{k=1}^{C n}    \widehat E  \left[ {\zeta_1 \over M_{n-r} }\right]\, \left[ 1- ( \widehat E(\Lambda_2(n-r)))^k \right], \end{eqnarray*}

\noindent where  the annealed expectation $\widehat E$ has the density $\zeta_2$ with respect to $ \q \otimes P $ and under $\widehat E$, $\Lambda_j$ are i.i.d and independent of $\zeta_1$.  Note that by the independence of $\zeta_2$ and $(\zeta_1, M_{n-r} )$ under $\q \otimes P $, we have $ \widehat E  \left[ {\zeta_1 \over M_{n-r} }\right] =  \q \otimes P   \left[       { \zeta_1  \over M_{n-r} } \right] \le C'$. 

To proceed, we employ the following
estimate, which will be proved below:
there exists a constant $c_1$ (that may depend on $\alpha$) so that
$\forall \ell \ge \ell_0$,
\begin{equation}\label{err2bis} 1- \widehat E( \Lambda_2(\ell))
\le e^{- c_1 \, \ell\,.}
\end{equation}
Since $n-r \ge \varepsilon n$, (\ref{err2bis}) yields that $I_{(\ref{err1bis})}  \to0$  as stated in  (\ref{err1bis}).

It remains to check (\ref{err2bis}). 
$\beta_\ell(o) - \beta(o)$ corresponds to the probability that
an excursion of the tree-valued walk  is higher  than $\ell$;
the latter is dominated by the probability that a level regeneration distance
is larger than $\ell$, which decays exponentially by \cite[Lemma 4.2(i)]{DGPZ}. It follows that 
 \begin{equation}
\label{campus}
    \p (\beta_\ell(o) - \beta(o) > e^{ - c_2 \ell}) \le e^{- c_2 \ell}\,, \qquad \forall \, \ell \ge \ell_0, 
\end{equation}
where $c_2$ may depend on $\alpha$.  Then $\e(  (\beta_\ell(o) - \beta(o)) \le 2 e^{- c_2 \ell}$.  Notice that for $R_1\le i < R_2$, $S_i <0$ hence $$\sum_{i=R_1}^{R_2 -1}    ( a_{S_i}^{(\ell)}- a_{S_i}^{(\infty)}) = \sum_{x\le 0}   ( a_{x}^{(\ell)}- a_{x}^{(\infty)}) (L_{R_2}^x- L_{R_1}^x) \le { 1\over \lambda}  \sum_{x\le 0}  \sum_{k=1}^{d_x^*-1} (\beta_\ell^{(x, k)}- \beta^{(x, k)}) (L_{R_2}^x- L_{R_1}^x),$$

 \noindent implying that $$ \q \otimes P  \left[ \sum_{i=R_1}^{R_2 -1}    ( a_{S_i}^{(\ell)}- a_{S_i}^{(\infty)})  \right] \le  { 2 \over \lambda} \q( d^*- 1) \, E  (R_2- R_1)   e^{- c_2 \ell}.  $$

By H\"{o}lder's inequality with ${1\over p}+{1\over q} =1$, 
\begin{eqnarray*}  \widehat E\left[1_{( \sum_{i=R_1}^{R_2 -1}    ( a_{S_i}^{(\ell)}- a_{S_i}^{(\infty)}) > e^{- c_3 \ell})} \right] &= &\q \otimes P  \left[\zeta_2 1_{( \sum_{i=R_1}^{R_2 -1}    ( a_{S_i}^{(\ell)}- a_{S_i}^{(\infty)}) > e^{- c_3 \ell})} \right]  
	\\ &\le&  ( \q \otimes P  ( (\zeta_2 )^p))^{1/p} \left(  \q \otimes P  (   \sum_{i=R_1}^{R_2 -1}    ( a_{S_i}^{(\ell)}- a_{S_i}^{(\infty)}) > e^{- c_3 \ell}) \right)^{1/q}  
	\\& \le& e^{- c_3 \ell} \,, \end{eqnarray*}  for some constant $c_3=c_3(\alpha, p, q, c_2)>0$ and for all large $\ell$. Now, using the  elementary inequality:  for any $j\ge 1$ and $x_1, ..., x_j \in [0, 1]$, $  1 - \prod_{i=1}^j (1-x_i) \le \sum_{i=1}^j x_i,$  we get    $$ 1-   \Lambda_2(\ell)  \le  \sum_{i=R_1}^{R_2 -1} { \lambda ( a_{S_i}^{(\ell)}- a_{S_i}^{(\infty)})   \over   1+ \lambda + \lambda a^{(\ell)}_{S_i}  } \le \sum_{i=R_1}^{R_2 -1}   ( a_{S_i}^{(\ell)}- a_{S_i}^{(\infty)})     . $$

\noindent Therefore $$  1- \widehat E( \Lambda_2(\ell)) \le 2\, e^{- c_3 \ell},  $$

\noindent proving  (\ref{err2bis}). The proof of the lemma is now complete. 
$\Box$

\bigskip

  We have the following representation of the velocity $v_\alpha$.
 \begin{theorem}[Velocity representation] \label{theo-3.4}
	 Assume \eqref{L:hyp1},  \eqref{L:hyp2}, \eqref{hyp-zeta2} and \eqref{L:hyp3}. Recall
	 the function  $f_\infty$, see \eqref{eq-fdef}. 
	 Then, $$    v_\alpha =  m   { \q \otimes P  \left[  \prod_{i=0}^{ R_1 -1}  f_{\infty}(S_i)   \,{\beta(u^*_1)\over M_\infty(u^*_1)}
\right]   \over  \q \otimes P  \left[  \prod_{i=0}^{ R_1 -1}  f_{\infty}(S_i)   \, {1\over M_\infty } 
\right]  }.  $$  
\end{theorem}

 We recall that $M_\infty \equiv M_\infty(u^*_0)$ and $u^*_0=o$.
 
{\noindent\bf Proof:}  Noticing that $\e( \Gamma_n(o))= m \, \e( \gamma_{n-1}(o))$. By \eqref{eq-gamma2}, \eqref{Gamma3} and Lemma \ref{lem-Udist}, we     immediately  obtain a representation of the velocity $v_\alpha$:    \begin{equation} \label{rep1}{ m\, \e(\beta) \over v_\alpha}=  { \q \otimes P  \left[  \prod_{i=0}^{ R_1 -1}  f_{\infty}(S_i)   \, {1\over M_\infty}
\right]   \over \q \otimes P    \left[  \prod_{i=R_1}^{ R_2 -1}  f_{\infty}(S_i)   \, \vert S_{R_2} - S_{R_1}\vert \right]} \, \, \sum_{y\ge 1} h(y) \,. \end{equation}

Going back to  (\ref{ebn1bis}), and  recalling that 
 $$ {m^r \over \lambda^r} \prod_{j=1}^r (1-  \beta_n(u^*_j))  = E  _{r, \omega} \Big( {\bf 1}_{( \tau_{S}(0) < \tau_{S}(n))}
    \prod_{i=0}^{\tau_S(0)-1}  f_n(S_i)\Big)   ,$$

\noindent 
we get that for any $r \le n-1$, $$ \e (B_n(o))= {m\over \lambda} \, \q \otimes E  _{r-1, \omega} \Big( {\bf 1}_{( \tau_{S}(0) < \tau_{S}(n))}
    \prod_{i=0}^{\tau_S(0)-1}  f_n(S_i) \, { \beta_n(u^*_r) \over M_n(u^*_r)}\Big) 
.$$

\noindent 
By shifting $P_{r-1, \omega} $ to $P_{0, \omega}$, we have that  for any $r \le n-1$, $$ \e (B_n(o))= {m\over \lambda} \, \q \otimes P\Big( {\bf 1}_{( \tau_{S}(1-r) < \tau_{S}(n-r+1))}
    \prod_{i=0}^{\tau_S(0)-1}  f_{n-r+1}(S_i) \, { \beta_{n-r+1}(u^*_1) \over M_{n-r+1}(u^*_1)}\Big) 
.$$ 

\noindent Repeating the renewal arguments in Section 3.2 which lead to Lemma \ref{lem-Udist} (the difference with $\Phi_n(r)$ only comes from the part before the regeneration time $R_1$), we see that  $$\lim_{n \to \infty } \e(B_n(o))= {m\over \lambda}\,    { \q \otimes P  \left[  \prod_{i=0}^{ R_1 -1}  f_{\infty}(S_i)   \, {\beta_\infty(u^*_1)\over M_\infty(u^*_1)}
\right]   \over \q \otimes P    \left[  \prod_{i=R_1}^{ R_2 -1}  f_{\infty}(S_i)   \, \vert S_{R_2} - S_{R_1}\vert \right]} \, \, \sum_{y\ge 1} h(y) \,. $$

\noindent On the other hand, $\lim_{n \to \infty } \e(B_n(o))= \e(B(o))= {m\over \lambda} \e(\beta(o))$. Comparing this with the velocity representation \eqref{rep1}, we get the result.  $\Box$

\medskip
Before applying Theorem \ref{theo-3.4}, we show that 
the conditions for the representation of $v_\alpha$ hold when $\alpha$ is small
enough. Recall our standing assumption that $p_0=0$, and the constant
$\kappa$, see \eqref{reg-exp}.

\begin{lemma}
	\label{lem-tildeqreg}
	There exists an $\alpha_0=\alpha_0(m,\kappa)$ such that
	if $0< \alpha<\alpha_0$ then \eqref{L:hyp1}, \eqref{L:hyp2}, \eqref{hyp-zeta2}, \eqref{unif} and \eqref{L:hyp3}
	hold.
\end{lemma}
\noindent
{\bf Proof:}
Note that $f_\infty(x)\leq (m^2+\lambda)/(m+m\lambda)$ and the right 
side is a bounded 
differentiable function of $\alpha$, which equals $1$ at 
$\alpha=0$. It follows that $f_\infty(x)\leq 1+c\alpha$ for
some constant $c=c(m)$. 

In what follows, we will make sure to use constants that do not depend
on $\alpha$.
Note that, since $\tilde P(\tau_S(1)=\infty)$ is bounded away
from $0$ uniformly in $\alpha$,
$$ \sum_{y=1}^\infty h(y)\leq 
C \sum_{n=0}^\infty  (1+c\alpha)^n \tilde P(R_1\geq n)
\leq C'\sum_{n=0}^\infty (1+c\alpha)^n e^{-\kappa n}\,,$$ 
where $C'=C'(\kappa,m)$ and we used \eqref{reg-exp}. In particular,
for $\alpha<\alpha_0(m,\kappa)$, we deduce \eqref{L:hyp1}.

The proof of \eqref{L:hyp2} is similar: since $|S_{R_2}-S_{R_1}|<R_2-R_1$,
the exponential moments \eqref{reg-exp} imply that it is enough to check
that $\q\otimes P [\zeta_1]<\infty$ and
$\q\otimes P [\zeta_2^{1+\delta}]<\infty$ 
for some $\delta>0$ independent of $\alpha$. Using again the estimate 
$f_\infty(x)\leq 1+c\alpha$ and the independence between $(M_k)_{1\le k \le \infty}$ and $(R_1, R_2)$, we see  that \eqref{L:hyp2}, \eqref{hyp-zeta2}  and \eqref{L:hyp3} 
follow  at once from \eqref{reg-exp} [for  \eqref{L:hyp3}, we also use the fact that $\q( { 1 \over M_{n-r}})=1, \forall n \ge r $].  It remains to check  the uniform integrabiltiy \eqref{unif}: Since $\zeta_1 \le (1+c \alpha)^{R_1}:= \zeta^*$, we have for any $a>0$ that \begin{eqnarray*}    \q \otimes P \left[   { \zeta_1 \over M_k} 1_{ ( { \zeta_1 \over M_k}  > a)} \right]  & \le   &   \q \otimes P \left[   { \zeta^* \over M_k} 1_{ (\zeta^* > a^{1/2})} \right] +   \q \otimes P \left[   { \zeta^* \over M_k} 1_{ ( { 1 \over M_k}  > a^{1/2})} \right] 
 	\\ & = & E \left[    \zeta^*   1_{ (\zeta^* > a^{1/2})} \right] +    E( \zeta^*) \, \q  \left[ { 1 \over M_k} 1_{ ( { 1 \over M_k}  > a^{1/2})} \right] ,\end{eqnarray*}

\noindent where we used the independence between $\zeta^*$ and $M_k
$, and $E$ denotes the expectation with respect to $P$.  Clearly, $E \left[    \zeta^*   1_{ (\zeta^* > a^{1/2})} \right]  =o(1)$ as $a \to \infty$. Observe that $ \q  \left[ { 1 \over M_k} 1_{ ( { 1 \over M_k}  > a^{1/2})} \right] = \p \left[  { 1 \over M_k}  > a^{1/2}  \right]  \le e^1 \, \e e^{ - a^{1/2} \, M_k} \le e^1 \, \e e^{ - a^{1/2}\, M_\infty}$. Since $p_0=0$,     $M_\infty >0, \p$-a.s.,  then $\e e^{ - a^{1/2}\, M_\infty} \to 0$ as $a \to \infty$, hence   $\q  \left[ { 1 \over M_k} 1_{ ( { 1 \over M_k}  > a^{1/2})} \right]  \to 0$ uniformly on  $k$ and we get \eqref{unif}.  \qed

\bigskip

{\noindent \bf Proof of Theorem \ref{theo-GW} (case $\alpha\searrow 0$).}
By proposition \ref{eresc}, 
$$
\lim_{\alpha\searrow 0} \frac{\e(\beta(x))}{\alpha}=\frac{{\cal
    D}^0}{2m}\,. 
$$

\noindent Then by the velocity representation for $v_\alpha$ in \eqref{rep1},  it is enough to prove that   \begin{equation}
		\label{new44} \lim_{\alpha \searrow 0}  { \q \otimes \widetilde
E \left[  \prod_{i=0}^{ R_1 -1}  f_{\infty}(S_i)   \, {1\over M_\infty}
\right]   \over \q \otimes P    \left[  \prod_{i=R_1}^{ R_2 -1}  f_{\infty}(S_i)   \, \vert S_{R_2} - S_{R_1}\vert \right]} \, \, \sum_{y\ge 1} h(y) =1\,.   \end{equation}

 Since we are interested in the limit $\alpha\searrow 0$, we may and
will assume throughout that $\alpha<\alpha_0(m,\kappa)$ the constant appeared in  Lemma \ref{lem-tildeqreg}.
We write in this proof $A\sim_\alpha B$ if 
$(A-B)/\alpha\to_{\alpha\searrow 0} 0$.

Note that  $f_\infty(x)\leq 1+c\alpha$ for some constant $c>0$.
Mimicking the proof of Lemma  \ref{lem-tildeqreg}, we therefore get that
$$ \q \otimes P    \left[  \prod_{i=R_1}^{ R_2 -1}  f_{\infty}(S_i)   
\, \vert S_{R_2} - S_{R_1}\vert \right] \sim_\alpha
E   \left[    \vert S_{R_2} - S_{R_1}\vert \right] = 
{ 1 \over P ( \tau_S(1)=\infty) }\sim_\alpha {m\over m-1}.$$

\noindent 
In the same way, $h(y) \sim_\alpha P ( \tau_S(-y) \le R_1 | 
\tau_S(1)=\infty)$ for $y \ge 1$, hence 
$$ \sum_{y\ge1}  h(y) \sim_\alpha 
E  ( |S_{R_1}|  | \tau_S(1)=\infty) \sim_\alpha {m\over m-1}.$$

\noindent 
Finally,  as $\alpha\to 0$,  $$\q \otimes P    \left[
\prod_{i=1}^{ R_1 -1}  f_{\infty}(S_i)   \, {1\over M_\infty} \right] \sim_\alpha
\q   \left[    {1\over M_\infty} \right] =1, $$

\noindent   implying \eqref{new44} and completing the proof of Theorem \ref{theo-GW}. $\Box$

 \bigskip
{\noindent\bf  Acknowledgment:} We thank Amir Dembo and Yuval Peres
 for useful discussions concerning regeneration times for Galton--Watson
 trees, Nina Gantert and Pierre Mathieu for discussing with some of us
 their paper \cite{GMP},
and Amir Dembo and Elie A\"{i}d\'{e}kon for comments
on an earlier version of this paper. We thank an anonymous referee for her/his comments.

\bigskip\bigskip


\medskip



\begin{thebibliography}{11}
	\bibitem{AB}
		D. J. Aldous and
		A.  Bandyopadhyay, 
		A survey of max-type recursive distributional equations,
		{\it Ann. Appl. Probab.} {\bf   15}  (2005),  pp.
                1047--1110.

\bibitem{DGPZ} A. Dembo, N. Gantert, Y. Peres and O. Zeitouni,
Large deviations for random walks on Galton--Watson trees:
averaging and uncertainty,
{\it Prob. Th. Rel. Fields} {\it 122} (2002), pp. 241--288. 

\bibitem{BD} A. Dembo, J.D. Deuschel, Markovian perturbation, response
  and fluctuation dissipation theorem, {\it Ann. I. H. Poincare,
    Prob. Stat.}, \textbf{46}, N. 3 (2010), 822-852. 

    \bibitem{DS} A. Dembo, N. Sun, Central limit theorem for biased 
	    random walk on multi-type Galton-Watson trees,
	    {\it arXiv:1011.4056v2}.

    \bibitem{Far} G. Faraud, A central limit theorem for random walk in random environment on marked Galton-Watson trees,
	    {\it arXiv:0812.1948v6}.
\bibitem{FHS} G. Faraud, Y. Hu and Z. Shi,
	Almost sure convergence for stochastically biased random walks on trees,
	{\it arXiv:1003:5505v2}.

\bibitem{Fel} W. Feller, {\it An introduction to probability theory and its
	applications}, 3rd Edition, John Wiley \& sons, 1968.

\bibitem{GMP} N. Gantert, P. Mathieu and  A. Piatnitski, 
Einstein relation for reversible diffusions in random environments,
{\it arXiv:1005:5665v2}.

\bibitem{OK} T. Komorowsky and S. Olla, 
Einstein relation for random walks in random environments,
{\it Stoch. Proc. Appl.} {\bf 115} (2005), pp. 1279--1301.

\bibitem{OK1} T. Komorowsky and S. Olla, 
On mobility and Einstein relation for tracers in time--mixing random
environments, {\it J. Stat. Phys} {\bf 118}, N. 3/4, (2005) pp. 407--435.

\bibitem{lero} J. L. Lebowitz and H. Rost,
 {The Einstein relation for the displacement of a test
 particle in a random environment},  {\it Stochastic Process. Appl.} \textbf{54}
    no. 2, (1994) 183--196.

\bibitem{LOU} M. Loulakis, {Einstein Relation for a tagged particle in
    simple exclusion processes.}
\emph{Comm.\ Math.\ Phys.}, {\bf 229}, (2002), pp. 347--367.

\bibitem{L97} R. Lyons,
A simple path to Biggins' martingale convergence for branching random walk. {\it Classical and modern branching processes (Minneapolis, MN, 1994)} 217--221,
IMA Vol. Math. Appl., 84, Springer, New York, 1997.

\bibitem{RL} R. Lyons, {Random walks and percolation on trees},
{\it Annals Probab.} {\bf 18} (1990), pp. 931--958.

\bibitem{LPP95}
R. Lyons, R. Pemantle and Y. Peres, ``Ergodic theory on
Galton--Watson trees:
speed of random walk
and dimension of harmonic measure",
{\it Ergodic Theory Dyn. Systems} {\bf 15} (1995), pp. 593--619.

\bibitem{LPP} R. Lyons, R. Pemantle and Y. Peres,
Biased random walks on Galton--Watson trees,
{\it Prob. Th. Rel. Fields} {\bf 106} (1996), pp. 249--264. 

\bibitem{PP}  R. Pemantle and Y. Peres,  The critical Ising model on trees, concave recursions and nonlinear capacity. {\it Ann. Probab. } {\bf 38 }  (2010),  184--206.

\bibitem{PZ} Y. Peres and O. Zeitouni,
A central limit theorem for biased random walks on Galton--Watson trees,
{\it Prob. Th. Rel. Fields} {\it 140} (2008), pp. 595--629.
\bibitem{stflour} O. Zeitouni, 
 {\it Random walks
in random environment}, XXXI Summer school in probability,
St Flour (2001). { Lecture notes in Math.} 1837 (Springer)
(2004), pp. 193--312.

\end{thebibliography}
\end{document}